\documentclass[11pt]{article}

\usepackage{times}
\usepackage[left=2.7cm, right=2.7cm, top=2.8cm, bottom=2.8cm]{geometry}
\usepackage[T1]{fontenc}
\usepackage[latin1]{inputenc}
\usepackage{amsmath,amssymb}
\usepackage{amsthm}
\usepackage{amscd}
\usepackage{mathdots}
\usepackage{multirow}
\usepackage[usenames,dvipsnames]{color}
\usepackage{epsfig}

\newtheorem{theorem}{Theorem}[section]

\newtheorem{corollary}[theorem]{Corollary}
\newtheorem{proposition}[theorem]{Proposition}

\theoremstyle{definition}
\newtheorem{example}[theorem]{Example}
\newtheorem{remark}[theorem]{Remark}
          
\numberwithin{table}{section}
\numberwithin{equation}{section}

\usepackage{color}

\begin{document}

\title{A Generalized Sampling Theorem for Stable Reconstructions in Arbitrary Bases} 

\author{%
{\sc
Ben Adcock}\\[2pt]
\small\em 
Department of Mathematics, Simon Fraser University, \\
\and
\sc Anders C. Hansen\\[2pt]
\small\em 
Department  of Applied Mathematics and Theoretical Physics, University of Cambridge.  \\
}

\date{}
\maketitle

\begin{abstract}
We introduce a generalized framework for sampling and reconstruction in separable Hilbert spaces.  Specifically, we establish that it is always possible to stably reconstruct a vector in an arbitrary Riesz basis from sufficiently many of its samples in any other Riesz basis.  This framework can be viewed as an extension of that of Eldar et al.  However, whilst the latter imposes stringent assumptions on the reconstruction basis, and may in practice be unstable, our framework allows for recovery in any (Riesz) basis in a manner that is completely stable.

Whilst  the classical Shannon Sampling Theorem is a special case of our theorem, this framework allows us to exploit additional information about the approximated vector (or, in this case, function), for example sparsity or regularity, to design a reconstruction basis that is better suited. Examples are presented illustrating this procedure.
\end{abstract}

 

\section{Introduction}
The Shannon Sampling Theorem, or the Nyquist--Shannon Sampling Theorem
as 
it is also called 
(we will refer to it as the NS-Sampling Theorem throughout the paper), 
is a mainstay in modern signal processing
and has become one of the most important theorems in mathematics of 
information \cite{unser2000sampling}. The list of applications of the theorem is long, and ranges from 
Magnetic Resonance Imaging (MRI) to 
sound engineering. We will in this paper address the question
on 
whether or not the NS-Sampling Theorem 
can be improved. In particular, given the same set of information,
could one design a reconstruction of a function that would be better
than that provided by the NS-Sampling Theorem. 
The answer to such a question will obviously depend on the 
type of functions considered.  However, suppose that we have some extra 
information about the functions to be reconstructed. One may, for
example, have information about a basis that is particularly suited for such functions.  
Could this information be used to improve the reconstruction given 
by the NS-Sampling Theorem, even if it is based on the same sampling 
procedure?  Although such a question has been posed before, and numerous extensions of the NS-Sampling Theorem have been developed \cite{aldroubi1996oblique,  Feichtinger1998,  Feichtinger2004,
Feichtinger2003, unser1994general}, the generalization we introduce in this paper is, to the best of our knowledge, a novel approach for this problem.

The well known  NS-Sampling Theorem \cite{jerri, book, Nyquist,Shannon, Whittaker} states that if
$$
f = \mathcal{F}g, \qquad g \in L²(\mathbb{R}),
$$ 
and $\mathrm{supp}(g) \subset [-T,T]$ for some $T > 0$, then both $f$
and $g$ can be reconstructed from point samples of $f$. In particular,
if $\epsilon \leq \frac{1}{2T}$ then 
$$
f(t) = \sum_{k=-\infty}^{\infty}f(k\epsilon)\mathrm{sinc}\left(\frac{t
    +k\epsilon}{\epsilon}\right) \quad L²\,\text{and unif. convergence,}
$$
$$
g(\cdot) = \epsilon \sum_{k=-\infty}^{\infty}f(k\epsilon)e^{2\pi \mathrm{i} \epsilon k \cdot }
\qquad L²\text{ convergence.}
$$ 
The quantity $\frac{1}{2T},$ which is the largest value of $\epsilon$
such that the theorem holds, is often referred to as the Nyquist
rate \cite{Nyquist}. In practice, when trying to reconstruct $f$ or $g$, one will most
likely not be able to access the infinite amount of information
required, namely, $\{f(k\epsilon)\}_{k \in \mathbb{Z}}.$  Moreover, even if we had access to all samples, we are limited by both processing power and storage to taking only a finite number. 
Thus, a more
realistic scenario is that one will be given a finite number of 
samples $\{f(k\epsilon)\}_{|k| \leq N}$, for some $N < \infty$, and seek to reconstruct $f$ from these samples.  The question is therefore: are the approximations
$$
f_N(\cdot)  = \sum_{k=-N}^{N}f(k\epsilon)\mathrm{sinc}\left(\frac{\cdot
    +k\epsilon}{\epsilon}\right), \qquad  g_N(\cdot) 
 = \epsilon \sum_{k=-N}^{N}f(k\epsilon)e^{2\pi \mathrm{i} \epsilon k \cdot} 
$$ 
optimal for $f$ and $g$ given the information 
$\{f(k\epsilon)\}_{|k| \leq N}$? To formalize this question consider 
the following.  For $N \in \mathbb{N}$ and $\epsilon > 0$, let
\begin{equation}\label{one}
\begin{split}
\Omega_{N,\epsilon} & = \{\xi \in \mathbb{C}^{2N+1}: \quad \xi =  \{f(k\epsilon)\}_{|k| \leq N},\, 
 f \in L^2(\mathbb{R})\cap C(\mathbb{R})\}.
\end{split}
\end{equation}
($C(\mathbb{R})$ denotes the set of continuous functions on $\mathbb{R}$).
Define the mappings (with a slight abuse of notation)
$$
\Lambda_{N,\epsilon,1}: \Omega_{N,\epsilon} 
\rightarrow L^2(\mathbb{R}), \quad  \Lambda_{N,\epsilon,2}: \Omega_{N,\epsilon} 
\rightarrow L^2(\mathbb{R}),
$$
\begin{equation}\label{two}
\begin{split}
 \Lambda_{N,\epsilon,1}(f) 
= \sum_{k=-N}^{N}f(k\epsilon)\mathrm{sinc}\left(\frac{\cdot
    +k\epsilon}{\epsilon}\right)
   \qquad 
\Lambda_{N,\epsilon,2}(f)
 = \epsilon \sum_{k=-N}^{N}f(k\epsilon)e^{2\pi \mathrm{i} \epsilon k \cdot}.
\end{split}
\end{equation}
The question is, 
given a class of functions $\Theta \subset L^2(\mathbb{R})$,
could there exist mappings 
$\Xi_{N,\epsilon,1}: \Omega_{N,\epsilon} \rightarrow L^2(\mathbb{R})$ and 
$\Xi_{N,\epsilon,2}: \Omega_{N,\epsilon} \rightarrow L^2(\mathbb{R})$ 
such that
\begin{equation*}
\begin{split}
\|\Xi_{N,\epsilon,1}(f) - f\|_{L^{\infty}(\mathbb{R})} < 
\|\Lambda_{N,\epsilon,1}(f) - f\|_{L^{\infty}(\mathbb{R})}
\quad \forall f, f = \mathcal{F}g, g \in \Theta,
\end{split}
\end{equation*} 
\begin{equation*}
\begin{split}
\|\Xi_{N,\epsilon,2}(f) - g\|_{L^{2}(\mathbb{R})} < 
\|\Lambda_{N,\epsilon,2}(f) - g\|_{L^{2}(\mathbb{R})} \quad \forall f, f = \mathcal{F}g, g \in \Theta.
\end{split}
\end{equation*} 
As we will see later, the answer to this question may very well be
yes, and the problem is therefore to find such mappings $\Xi_{N,\epsilon,1}$ and
$\Xi_{N,\epsilon,2}.$ 

\begin{figure}
\centering
\includegraphics[height=48mm]{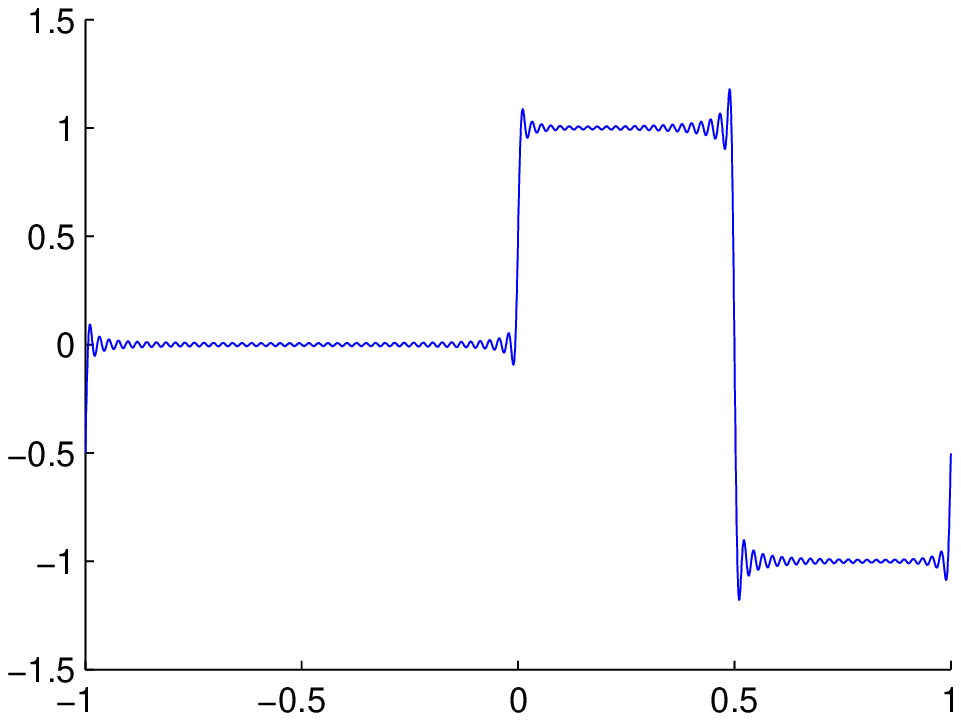}
\includegraphics[height=48mm]{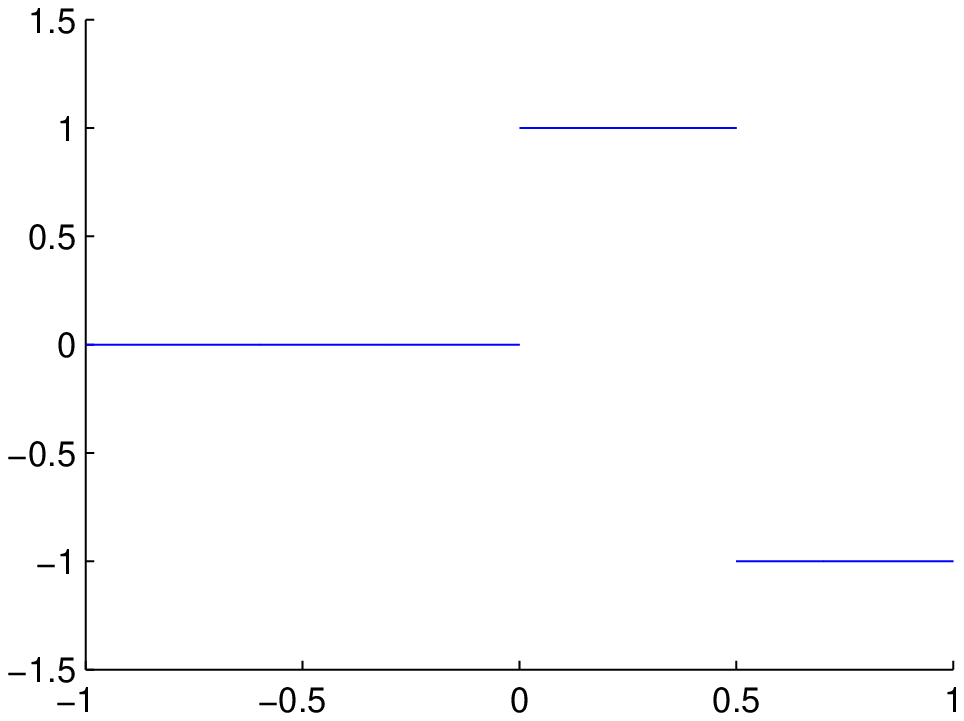}
\caption{The figure shows $\Lambda_{N,\epsilon,2}(f)$ for $f = \mathcal{F}g$, $N = 500$ and 
$\epsilon = 0.5$ (left) as well as $g$ (right).}
\label{freddy}
\end{figure}

As motivation for our work, consider the following
 reconstruction problem. 
Let $g$ be defined by 
$$
g(t) = 
\begin{cases} 1 & t \in [0,1/2)\\
-1&  t \in [1/2,1] \\
0 & t \in \mathbb{R}\setminus [0,1].
\end{cases}
$$ 
This is the well-known Haar wavelet.  Due to the discontinuity, there is 
no way one can exactly reconstruct this function with only finitely many
function samples if one insists on using the mapping
$\Lambda_{N,\epsilon,2}$. We have visualized the reconstruction of $g$
using $\Lambda_{N,\epsilon,2}$ in Figure \ref{freddy}.  In addition to $g$ not being reconstructed exactly, the approximation $\Lambda_{N,\epsilon,2}(g)$ is polluted by oscillations near the discontinuities of $g$.  Such oscillations are indicative of the well-known Gibbs phenomenon in recovering discontinuous signals from samples of their Fourier transforms \cite{jerrigibbs}.  This phenomenon is a major hurdle in many applications, including image and signal processing.  Its resolution has, and continues to be, the subject of significant inquiry \cite{tadmoracta}.

However, it is tempting to think that one
could construct a mapping $\Xi_{N,\epsilon,2}$ that would yield a
better result. Suppose for a moment that we do not know $g$, but we do
have some extra information. In particular, suppose that we know that 
$g \in \Theta$, where 
\begin{equation}\label{three}
\Theta = \left\{h \in L^2(\mathbb{R}): h = \sum_{k = 1}^M \beta_k
  \psi_k\right\}, 
\end{equation} 
for some finite number $M$ and where $\{\psi_k\}$ are the Haar wavelets on the interval $[0,1].$ 
Could we, based on the extra knowledge of $\Theta$, construct 
mappings $\Xi_{N,\epsilon,1}: \Omega_{N,\epsilon} \rightarrow
L^2(\mathbb{R})$  and $\Xi_{N,\epsilon,2}: \Omega_{N,\epsilon} \rightarrow
L^2(\mathbb{R})$ such that 
\begin{equation*}
\begin{split}
\sup\{\|\Xi_{N,\epsilon,1}(f) - f\|_{L^{\infty}(\mathbb{R})}: g \in \Theta ,
f=\mathcal{F}g\} &< 
\sup\{\|\Lambda_{N,\epsilon,1}(f) - f\|_{L^{\infty}(\mathbb{R})}: g \in \Theta, f=\mathcal{F}g \},\\
\sup\{\|\Xi_{N,\epsilon,2}(f) - g\|_{L^{2}(\mathbb{R})}: g \in \Theta ,
f=\mathcal{F}g\} &< 
\sup\{\|\Lambda_{N,\epsilon,2}(f) - g\|_{L^{2}(\mathbb{R})}: g \in \Theta, f=\mathcal{F}g \}?
\end{split}
\end{equation*} 
Indeed, this is the case, and a consequence of our framework is that 
it is possible to find $\Xi_{N,\epsilon,1}$ and $\Xi_{N,\epsilon,2}$ such
that
\begin{equation*}
\begin{split}
\sup\{\| \Xi_{N,\epsilon,1}(f) - f\|_{L^{\infty}(\mathbb{R})}: g \in \Theta ,
f=\mathcal{F}g\} &= 0,\\
\sup\{\| \Xi_{N,\epsilon,2}(f) - g\|_{L^{2}(\mathbb{R})}: g \in \Theta ,
f=\mathcal{F}g\} &= 0,
\end{split}
\end{equation*}
provided $N$ is sufficiently large. In other words, one gets perfect reconstruction. Moreover, the reconstruction is done in a completely stable way.  

The main tool for this task is a generalization of
the NS-Sampling Theorem that allows reconstructions in arbitrary
bases.  Having said this, whilst the Shannon Sampling Theorem is the most standard example, the framework we develop addresses the more abstract problem of recovering a vector (belonging to some separable Hilbert space $\mathcal{H}$) given a finite number of its samples with respect any Riesz basis of $\mathcal{H}$.

\subsection{Organization of the Paper}
We have organized the paper as follows.  In Section \ref{notation} we introduce notation and idea of finite sections 
of infinite matrices, a concept that will be crucial throughout the paper, and in Section \ref{early_work} we discuss existing literature on this topic, including the work of Eldar et al \cite{eldar2003sampling,eldar2005general,unser1994general}. The main theorem is presented and proved in Section \ref{new}, where we also show the connection to the classical NS-Sampling Theorem. The error bounds in the generalized sampling theorem involve several important constants, which can be estimated numerically. We therefore devote Section \ref{Normbounds} to discussions on how to compute crucial constants and functions that are useful for providing error estimates. Finally, in Section \ref{Numerics} we provide several examples to support the generalized sampling theorem and to justify our approach. 

\section{Background and Notation}\label{notation}
Let $\mathrm{i}$ denote the imaginary unit.  Define the Fourier transform $\mathcal{F}$ by 
$$
(\mathcal{F}f)(y) = \int_{\mathbb{R}^d} f(x) e^{-2\pi \mathrm{i} x \cdot y}\, dx,
\qquad f \in L^1(\mathbb{R}^d), 
$$
where, for vectors 
$x, y \in \mathbb{R}^d,$ $x \cdot y = x_1y_1 + \hdots + x_dy_d.$ 
Aside from the Hilbert space $L^2(\mathbb{R}^d)$, we now introduce
two other important Hilbert spaces: namely, 
$$
l^2(\mathbb{N}) = \left\{\alpha = \{\alpha_1, \alpha_2, \hdots\}: 
\sum_{k \in \mathbb{N}} |\alpha^2_k| < \infty \right\}
$$
and 
$$
l^2(\mathbb{Z}) = \left\{\beta = 
\{\hdots \beta_{-1}, \beta_0, \beta_{1} \hdots\}: 
\sum_{k \in \mathbb{Z}} |\beta^2_k| < \infty \right\},
$$
with their obvious inner products. We will also consider abstract Hilbert spaces. In this case we will use the 
notation $\mathcal{H}$. Note that $\{e_j\}_{j \in \mathbb{N}}$ and 
$\{e_j\}_{j \in \mathbb{Z}}$ will always denote the natural
bases for $l^2(\mathbb{N})$ and $l^2(\mathbb{Z})$ respectively. We may also use the notation 
$\mathcal{H}$ for both $l^2(\mathbb{N})$ and $l^2(\mathbb{Z})$ (the meaning will be clear from the context).  Throughout the paper, the symbol $\otimes$ will denote the standard tensor product on Hilbert spaces.

The concept of infinite matrices will be quite crucial in the theory, 
and also finite sections of such matrices. We will consider infinite
matrices as operators from both $l^2(\mathbb{N})$ to $l^2(\mathbb{Z})$ and 
$l^2(\mathbb{N})$ to $l^2(\mathbb{N})$. The set of bounded
operators from a Hilbert space $\mathcal{H}_1$ to a Hilbert space 
$\mathcal{H}_2$ will be denoted by
$\mathcal{B}(\mathcal{H}_1,\mathcal{H}_2)$. As infinite matrices are 
unsuitable for computations we must reduce any infinite
matrix to a more tractable finite-dimensional object. The standard means in which to do this is
via finite sections. In particular, let 
\begin{equation*}
U =
\left(
\begin{matrix}
\vdots            & \vdots     & \vdots   &  \iddots\\
 u_{-1,1}         & u_{-1,2}    & u_{-1,3}  & \hdots \\
 u_{0,1}         & u_{0,2}    & u_{0,3}  & \hdots \\
 u_{1,1}         & u_{1,2}    & u_{1,3}  & \hdots \\
\vdots          & \vdots    & \vdots  & \ddots\\
\end{matrix}
\right), \quad U \in \mathcal{B}(l^2(\mathbb{N}),l^2(\mathbb{Z})).
\end{equation*}
For $n \in \mathbb{N}$, define $P_n$ to be the
projection onto $\mathrm{span}\{e_1,\hdots,e_n\}$ and, for odd 
$m \in \mathbb{N}$, let $\widetilde P_m$ be the projection onto 
$\mathrm{span}\{e_{-\frac{m-1}{2}},\hdots, e_{\frac{m-1}{2}}\}$.
Then $\widetilde P_m U P_n$ may be interpreted as 
$$
\left(
\begin{matrix}
 u_{-\frac{m-1}{2},1}         & \hdots    & u_{-\frac{m-1}{2},n}\\
\vdots            & \vdots     & \vdots  \\

 u_{\frac{m-1}{2},1}         & \hdots    & u_{\frac{m-1}{2},n}\\
\end{matrix}
\right),
$$
a $m \times n$ finite section of $U$. 
 Finally, the spectrum of any operator $T \in \mathcal{B}(\mathcal{H})$ will 
be denoted by $\sigma(T).$

\section{Connection to Earlier Work}\label{early_work} 
The idea of reconstructing signals in arbitrary bases is certainly not a new idea and this topic has gone through extensive investigations in the last several decades.  The papers by Unser and Aldroubi \cite{aldroubi1996oblique,unser1994general} have 
been very influential and these ideas have been generalized to arbitrary Hilbert spaces by Eldar \cite{eldar2003sampling, eldar2005general}. The 
abstract framework introduced by Eldar is very powerful because of its general nature. Our framework is 
based on similar generalizations, yet it incorporates several key distinctions, resulting in a number of advantages.
Before introducing this framework, let us first review  some of the key concepts of \cite{eldar2005general}.

Let $\mathcal{H}$ be a separable Hilbert space and let $f \in \mathcal{H}$ be an element we would like to reconstruct from some measurements. Suppose that we are given linearly independent 
sampling vectors $\{s_k\}_{k\in\mathbb{N}}$ that span a subspace 
$\mathcal{S} \subset \mathcal{H}$ and form a Riesz basis, and assume that we can access the sampled inner products 
$c_k = \langle s_k, f \rangle$, $k=1,2 \hdots$. Suppose also that we are given linearly independent 
reconstruction vectors 
$\{w_k\}_{k\in\mathbb{N}}$ that span a subspace $\mathcal{W} \subset \mathcal{H}$ and also form a Riesz basis. The task is to obtain
a reconstruction $\tilde f \in \mathcal{W}$ based on the sampling data $\{c_k\}_{k\in\mathbb{N}}$. The natural choice, as suggested 
in \cite{eldar2005general}, is 
\begin{equation}\label{el}
\tilde f  = W(S^*W)^{-1}S^*f,
\end{equation} 
where the so-called \textit{synthesis operators} $S, W:l^2(\mathbb{N}) \rightarrow \mathcal{H}$ are defined by
\begin{equation*}
Sx = x_1s_1 + x_2s_2 + \hdots, \qquad Wy = y_1w_1 +y_2 w_2 +\hdots ,
\end{equation*}
and their adjoints $S^*, W^*:\mathcal{H} \rightarrow l^2(\mathbb{N})$ are easily seen to be
$$
S^*g = \{\langle s_1,g\rangle,  \langle s_2,g\rangle, \hdots \}, \qquad W^*h = \{\langle w_1,h\rangle, \langle w_2,h\rangle \hdots \}.
$$
Note that $S^*W$ will be invertible if and only if 
$$
\mathcal{H} = \mathcal{W} \oplus \mathcal{S}^{\perp}.
$$
This gives a very nice and intuitive abstract formulation of the reconstruction.  However, in practice we will never have the luxury of being able to acquire nor process the infinite amount of samples  $\langle s_k, f \rangle$, $k=1,2 \hdots$, needed to construct $\tilde{f}$. An important question to ask is therefore:
\begin{itemize}
\item[] What if we are given only the first $m \in \mathbb{N}$ samples
$\langle s_k, f \rangle$, $k=1, \hdots, m$?  In this case we cannot use (\ref{el}), and we may ask: what to do?
\end{itemize}
Fortunately, there is a simple finite dimensional analog to the infinite dimensional ideas discussed above.
Suppose that we are given $m \in \mathbb{N}$ linearly independent 
sampling vectors $\{s_1,\hdots, s_m\}$ that span a subspace 
$\mathcal{S}_m \subset \mathcal{H}$, and assume that we can access the sampled inner products 
$c_k = \langle s_k, f \rangle$, $k=1,\hdots,m$. Suppose also that we are given linearly independent 
reconstruction vectors 
$\{w_1, \hdots, w_m\}$ that span a subspace $\mathcal{W}_m \subset \mathcal{H}$. The task is to construct an 
approximation $\tilde f \in \mathcal{W}_m$ to $f$ based on the samples $\{c_k\}_{k=1}^m.$ In particular, we 
are interested in finding coefficients $\{d_k\}_{k=1}^m$ (that are computed from the samples 
$\{c_k\}_{k=1}^m$) such that $\tilde f = \sum_{k=1}^m d_k w_k$. The reconstruction suggested in \cite{eldar2003FAA} is 
\begin{equation}\label{f}
\tilde f = \sum_{k=1}^m d_k w_k = W_m(S_m^*W_m)^{-1}S_m^*f,
\end{equation}
where the operators $S_m, W_m : \mathbb{C}^m \rightarrow \mathcal{H}$ are defined by 
\begin{equation}\label{S_and_W}
S_mx = x_1s_1 + \hdots  +x_ms_m, \qquad W_my = y_1w_1 + \hdots +y_m w_m,
\end{equation}
and their adjoints $S^*, W^*:\mathcal{H} \rightarrow \mathbb{C}^m$ are easily seen to be
$$
S_m^*g = \{\langle s_1,g\rangle, \hdots , \langle s_m,g\rangle\}, \qquad W_m^*h = \{\langle w_1,h\rangle, \hdots , 
\langle w_m,h\rangle\}.
$$
From this it is clear that we can express $S_m^*W_m: \mathbb{C}^m \rightarrow \mathbb{C}^m$ as 
the matrix
\begin{equation}\label{the_matrix}
\left(
\begin{matrix}
\langle s_1,w_1 \rangle         & \hdots  &  \langle s_1,w_m \rangle \\
\vdots            & \vdots     & \vdots  \\
\langle s_m,w_1 \rangle          & \hdots    & \langle s_m ,w_m\rangle \\
\end{matrix}
\right).
\end{equation}
Also, $S_m^*W_m$ is invertible if and only if  and (\cite[Prop. 3]{eldar2003FAA})
\begin{equation}\label{SW}
\mathcal{W}_m \cap \mathcal{S}_m^{\perp} = \{0\}.
\end{equation} 
Thus, to construct $\tilde f$ one simply solves a linear system of equations.
The error can now conveniently be bounded from above and below by 
$$
\|f - P_{\mathcal{W}_m}f\| \leq \|f - \tilde f\| \leq \frac{1}{\cos(\theta_{\mathcal{W}_m\mathcal{S}_m})}\|f - P_{\mathcal{W}_m}f\|,
$$ 
where $P_{\mathcal{W}_m}$ is the projection onto $\mathcal{W}_m$,
$$
\cos(\theta_{\mathcal{W}_m\mathcal{S}_m}) = \mathrm{inf}\{\|P_{\mathcal{S}_m}g\|: g \in \mathcal{W}_m, \|g\| = 1\}
$$
is the cosine of the angles between the subspaces $\mathcal{S}_m$ and $\mathcal{W}_m$
and $P_{\mathcal{S}_m}$ is the projection onto $\mathcal{S}_m$ \cite{eldar2003FAA}.  Note that if $f \in \mathcal{W}_m$, then $\tilde{f} = f $ exactly, a feature known as \textit{perfect} recovery.  Another facet of this framework is so-called \textit{consistency}: the samples $\langle s_{j} ,\tilde{f}\rangle$, $j=1,\ldots,m$, of the approximation $\tilde{f}$ are identical to those of the original function $f$ (indeed, $\tilde{f}$, as given by (\ref{f}), can be equivalently defined as the unique element in $\mathcal{W}_{m}$ that is consistent with $f$).

Returning to this issue at hand, there are now several important questions to ask:
\begin{itemize}
\item[(i)] What if $\mathcal{W}_m \cap \mathcal{S}_m^{\perp} \neq \{0\}$ so that $S_m^*W_m$ is not invertible? It is very easy to 
construct theoretical examples such that $S_m^*W_m$ is not invertible, however (as we will see below), such situations may very well occur in applications. In fact, $\mathcal{W}_m \cap \mathcal{S}_m^{\perp} = \{0\}$ is a rather strict condition. If we have that $\mathcal{W}_m \cap \mathcal{S}_m^{\perp} \neq \{0\}$ does that mean that is is impossible to construct an approximation $\tilde f$ from the samples $S_m^*f$?
\item[(ii)] What if $\|(S_m^*W_m)^{-1}\|$ is large? The stability of the method must clearly depend on the quantity 
$\|(S_m^*W_m)^{-1}\|$. Thus, even if $(S_m^*W_m)^{-1}$ exists, one may not be able to use the method in practice as there will likely be increased sensitivity to both round-off error and noise.
\end{itemize}
Our framework is specifically designed to tackle these issues.
But before we present our idea, let us consider some examples where the issues in (i) and (ii) will be present.
\begin{figure}
\centering
\includegraphics[height=39mm]{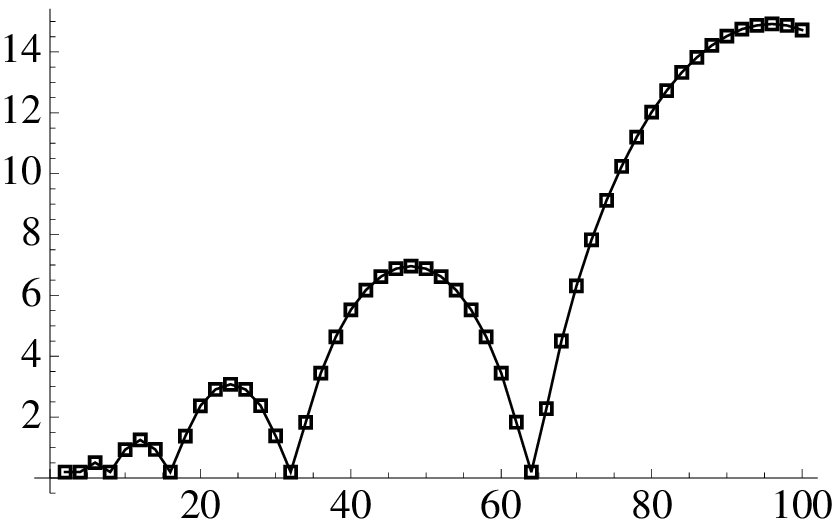} \hspace{2pc}
\includegraphics[height=39mm]{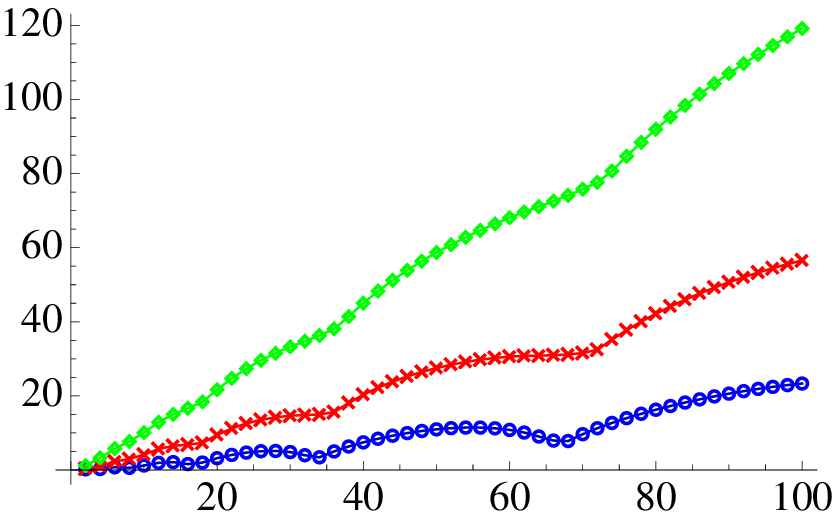}
\caption{This figure shows $\log_{10}\|(S^*_{\epsilon,m}W_m)^{-1}\|$ as a function of $m$ and $\epsilon$ 
for $m = 1,2,\ldots,100$. The left plot corresponds to $\epsilon = 1$, whereas the right plot corresponds to 
$\epsilon = 7/8$ (circles), $\epsilon = 1/2$ (crosses) and $\epsilon = 1/8$ (diamonds).}
\label{freddy33}
\end{figure}

\begin{example}
As for (i), the simplest example is to let $\mathcal{H} = l^2(\mathbb{Z})$ and $\{e_j\}_{j\in \mathbb{Z}}$ 
be the natural basis ($e_j$ is the infinite sequence with $1$ in its $j$-th coordinate and zeros elsewhere). 
For $m \in \mathbb{N}$, let the sampling vectors $\{s_k\}_{k=-m}^m$ and the reconstruction vectors 
$\{w_k\}_{k=-m}^m$ be defined by $s_k = e_k$ and $w_k = e_{k+1}.$ Then, clearly, $\mathcal{W}_m \cap \mathcal{S}_m^{\perp} =\mathrm{span}\{e_{m+1}\}$.  
\end{example}

\begin{example}\label{exp}
For an example of more practical interest, consider the following:  For $0 < \epsilon \leq 1$ let
$\mathcal{H} = L^2([0,1/\epsilon])$, and, for 
odd $m \in \mathbb{N}$, define 
the sampling vectors 
$$
\{s_{\epsilon,k}\}_{k=-(m-1)/2}^{(m-1)/2}, \qquad s_{\epsilon,k} = e^{-2\pi \mathrm{i}\epsilon k \cdot}\chi_{[0,1/\epsilon]},
$$ 
(this is exactly the type of measurement 
vector that will be used if one models Magnetic Resonance Imaging) and 
let the reconstruction vectors $\{w_k\}_{k=1}^m$ denote the $m$ first Haar wavelets on 
$[0,1]$ (including the constant function, $w_1 = \chi_{[0,1]}$). Let
 $S_{\epsilon,m}$ and $W_m$ be as in (\ref{S_and_W}), 
 according to the sampling and reconstruction vectors just defined.
A plot of $\|(S_{\epsilon,m}^*W_m)^{-1}\|$ as a function of $m$ and $\epsilon$ is given in Figure \ref{freddy33}.  As we observe, for $\epsilon =1$ only certain values of $m$ yield stable reconstruction, whereas for the other values 
of $\epsilon$ the quantity $\|(S_{\epsilon,m}^*W_m)^{-1}\|$ grows exponentially with $m$, making the problem severely ill-conditioned. Further computations suggest that $\|(S_{\epsilon,m}^*W_m)^{-1}\|$ increases exponentially with $m$ not just for these values of $\epsilon$, but for all $0<\epsilon<1.$ 
\end{example}

\begin{example}\label{exp2}
Another example can be made by replacing the Haar wavelet basis with the basis consisting of Legendre polynomials (orthogonal polynomials on $[-1,1]$ with respect to the Euclidean inner product).

In Figure \ref{Legfig} we plot the quantity $\| (S_{\epsilon,m}^* W_m)^{-1} \|$.  Unlike the previous example, this quantity grows exponentially and monotonically in $m$.  Whilst this not only makes the method highly susceptible to round-off error and noise, it can also prevent convergence of the approximation $\tilde{f}$ 
(as $m \rightarrow \infty$).  In essence, for convergence to occur, the error $\| f - P_{\mathcal{W}_m} f \|$ must decay more rapidly than the quantity $\| (S_{\epsilon,m}^* W_m)^{-1} \|$ grows.  Whenever this is not the case, convergence is not assured.  To illustrate this shortcoming, in Figure \ref{Legfig} we also plot the error $\| f - \tilde{f} \|$, where $f(x)  =\frac{1}{1+16 x^2}$.  The complex singularity at $x = \pm \frac{1}{4} \mathrm{i}$ limits the convergence rate of $\| f - P_{\mathcal{W}_m} f \|$ sufficiently so that $\tilde{f}$ does not converge to $f$.  Note that this effect is well documented as occurring in a related reconstruction problem, where a function defined on $[-1,1]$ is interpolated at $m$ equidistant pointwise samples  by a polynomial of degree $m-1$.  This is the famous Runge phenomenon.  The problem considered above (reconstruction from $m$ Fourier samples) can be viewed as a continuous analogue of this phenomenon.

\begin{figure}
\centering
\includegraphics[height=39mm]{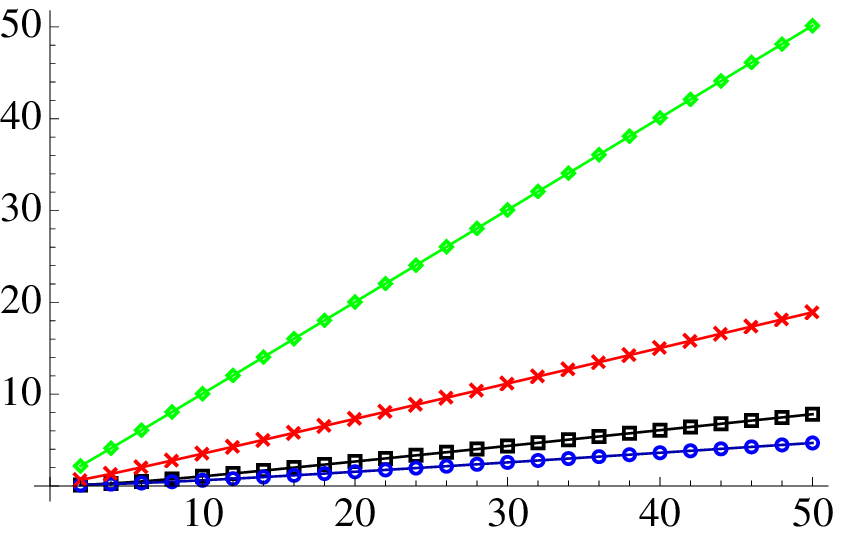}  \hspace{2pc}
 \includegraphics[height=39mm]{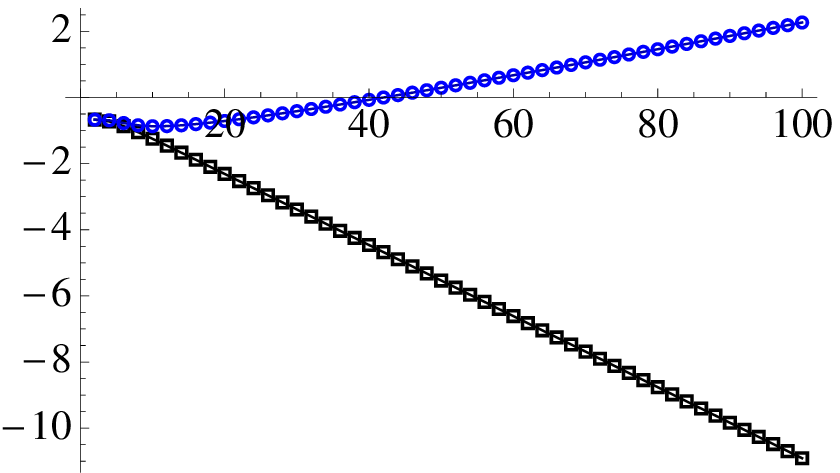}
\caption{The left figure shows $\log_{10}\|(S_{\epsilon,m}^*W_m)^{-1}\|$ as a function of $m$ for $m=2,4,\ldots,50$ and $\epsilon = 1,\frac{7}{8},\frac{1}{2},\frac{1}{8}$ (squares, circles, crosses and diamonds respectively).  The right figure shows $\log_{10}\| f - P_{\mathcal{W}_{m}} f \|$ (squares) and $\log_{10}\| f - \tilde{f} \|$ (circles) for $m= 2,4,6,\ldots,100$, where $f(x) = \frac{1}{1+16 x^2}$.}
\label{Legfig}
\end{figure}
\end{example}

Actually, the phenomenon illustrated in Examples \ref{exp} and \ref{exp2} is not hard to explain if one looks at the problem from an operator-theoretical point of view. This is the topic of the next section. 

\subsection{Connections to the Finite Section Method}
To illustrate the idea, let  $\{s_k\}_{k\in \mathbb{N}}$ and $\{w_k\}_{k\in \mathbb{N}}$ be two sequences 
of linearly independent elements in a Hilbert space $\mathcal{H}$. Define the infinite matrix $U$ by 
\begin{equation}\label{inf_U}
U =
\left(
\begin{matrix}
 u_{11}         & u_{12}    & u_{13}  & \hdots \\
 u_{21}         & u_{22}    & u_{23}  & \hdots \\
 u_{31}         & u_{32}    & u_{33}  & \hdots \\
\vdots          & \vdots    & \vdots  & \ddots\\
\end{matrix}
\right), \qquad u_{ij} =   \langle s_i,w_j\rangle.
\end{equation}
Thus, by (\ref{the_matrix}) the operator $S_m^*W_m$ is simply the $m$ by $m$ finite section of 
$U$. In particular 
$$
S_m^*W_m = P_mUP_m\vert_{P_m l^2(\mathbb{N})}.
$$ 
The finite section method has been studied extensively for the last decades 
\cite{bottcher1996, Silbermann2001, hansen2008, Lindner2006}. 
It is well known that even if 
$U$ is invertible then $P_mUP_m\vert_{P_m l^2(\mathbb{N})}$ may never be invertible for any $m$. 
In fact one must have rather strict conditions on $U$ for $P_mUP_m\vert_{P_m l^2(\mathbb{N})}$ 
to be invertible with uniformly bounded inverse (such as positive self-adjointness, for example \cite{Lindner2006}). 
In addition, even if 
$U:l^2(\mathbb{N}) \rightarrow l^2(\mathbb{N})$ is invertible and 
$P_mUP_m\vert_{P_m l^2(\mathbb{N})}$ is invertible for all $m \in \mathbb{N},$ it may be the case that, if
$$
x=U^{-1}y,  \quad x,y \in l^2(\mathbb{N}), \qquad x_m =  (P_mUP_m\vert_{P_m l^2(\mathbb{N})})^{-1}P_my
$$
then
$$
x_m \nrightarrow x, \qquad m \rightarrow \infty.
$$

Suppose that 
$\{s_k\}_{k\in \mathbb{N}}$ and $\{w_k\}_{k\in \mathbb{N}}$ are two Riesz bases for closed subspaces 
$\mathcal{S}$ and $\mathcal{W}$ of a separable Hilbert space $\mathcal{H}$. 
Define the operators $S, W:l^2(\mathbb{N}) \rightarrow \mathcal{H}$ by
\begin{equation}\label{SandW}
Sx = x_1s_1 +x_2s_2 + \hdots , \qquad Wy = y_1w_1 + y_2w_2+ \hdots.
\end{equation}
Suppose now that 
$(S^*W)^{-1}$ exists.  For $m \in \mathbb{N}$, let the spaces 
$\mathcal{S}_m, \mathcal{W}_m$ and operators $S_m, W_m : \mathbb{C}^m \rightarrow \mathcal{H}$ be defined 
as in Section \ref{early_work} according to the vectors $\{s_k\}_{k=1}^m$ and $\{w_k\}_{k=1}^m$ respectively. 
The following scenarios may now arise:
\begin{itemize}
\item[(i)] $\mathcal{W} \cap \mathcal{S}^{\perp} = \{0\}$, yet
$$
\mathcal{W}_m \cap \mathcal{S}^{\perp}_m \neq \{0\}, \qquad \forall \, m \in \mathbb{N}.
$$
\item[(ii)] $\|(S^*W)^{-1}\| < \infty$ and the inverse $(S_m^*W_m)^{-1}$ exists for all $m \in \mathbb{N}$, but
$$
\|(S_m^*W_m)^{-1}\| \longrightarrow \infty, \qquad m \rightarrow \infty.
$$ 
\item[(iii)] $(S_m^*W_m)^{-1}$ exists for all $m \in \mathbb{N},$ however
$$
W_m(S_m^*W_m)^{-1}S_m^*f \nrightarrow f, \qquad  m \rightarrow \infty,
$$
for some $f \in \mathcal{W} $.
\end{itemize}

Thus, in order for us to have a completely general sampling theorem we must try to
extend the framework described in this section in order to overcome the obstacles listed above.

\section{The New Approach}\label{new}
\subsection{The Idea}
One would like to have a completely general sampling theory that can be described as follows: 
\begin{itemize}
\item[(i)]
We have a signal 
$f \in \mathcal{H}$ and a Riez basis $\{w_k\}_{k\in\mathbb{N}}$ that spans some closed subspace 
$\mathcal{W} \subset \mathcal{H}$, and 
$$
f = \sum_{k=1}^{\infty} \beta_k w_k, \qquad \beta_k \in \mathbb{C}.
$$ 
So $f \in \mathcal{W}$ (we may also typically have some information on the decay rate of the $\beta_k$s, 
however, this is not crucial for our theory).
\item[(ii)]
We have sampling vectors $\{s_k\}_{k\in\mathbb{N}}$ that form a Riez basis for a closed subspace 
$\mathcal{S} \subset \mathcal{H},$ (note that 
we may not 
have the luxury of choosing such sampling vectors as they may be specified by some particular model, as is 
the case in MRI) and we can access the sampling values $\{\langle s_k,f\rangle\}_{k\in \mathbb{N}}.$
\end{itemize}
{\bf Goal:} reconstruct the best possible approximation $\tilde f \in \mathcal{W}$ based on the finite subset $\{\langle s_k,f\rangle\}_{k=1}^m$ of the sampling information  $\{\langle s_k,f\rangle\}_{k\in \mathbb{N}}$.  

We could have chosen $m$ vectors $\{w_1, \hdots, w_m\}$ and defined the operators $S_m$ and $W_m$ as in 
(\ref{S_and_W}) (from $\{w_1, \hdots, w_m\}$ and $\{s_1, \hdots, s_m\}$) and let $\tilde f$ be defined by (\ref{f}). However, this may be impossible as $S_m^*W_m$ may not be invertible (or the inverse may have a very large norm), as discussed in Examples 
\ref{exp} and \ref{exp2}. The question is then: what to do? 

To deal with these issues we will launch an abstract sampling theorem that extends the ideas discussed above.
To do so, we first notice that, since $\{ s_j \}$ and $\{ w_j \}$ are Riesz bases, there exist constants $A,B,C,D > 0$ such that 
\begin{equation}\label{riesz}
\begin{split}
&A \sum_{k\in\mathbb{N}}|\alpha_k|^2 \leq \left\| \sum_{k\in\mathbb{N}} \alpha_kw_k \right\|^2 \leq 
B \sum_{k\in\mathbb{N}}|\alpha_k|^2\\
&C \sum_{k\in\mathbb{N}}|\alpha_k|^2 \leq \left\| \sum_{k\in\mathbb{N}} \alpha_k s_k \right\|^2 \leq 
D \sum_{k\in\mathbb{N}}|\alpha_k|^2, \quad \forall \,\{\alpha_1, \alpha_2, \hdots\} \in l^2(\mathbb{N}).
\end{split}
\end{equation}
Now let $U$ be defined as in (\ref{inf_U}).
Instead of dealing with $P_mUP_m\vert_{P_m l^2(\mathbb{N})} = S_m^*W_m$ we propose to choose $n \in \mathbb{N}$ and compute the solution 
$\{\tilde \beta_1,\hdots, \tilde \beta_n\}$ of the following equation:
\begin{equation}\label{the_eq}
A\left(
\begin{matrix}
 \tilde \beta_1 \\
 \tilde \beta_2 \\
\vdots       \\
\tilde \beta_n\\
\end{matrix}
\right)
=
P_nU^*P_m
\left(
\begin{matrix} 
\langle s_1,f\rangle  \\
 \langle s_2,f\rangle \\
 \vdots  \\
\langle s_m,f \rangle  \\
 \end{matrix}
\right), \qquad A = P_nU^*P_mUP_n\lvert_{P_n\mathcal{H}},
\end{equation}
provided a solution exists (later we will provide estimates on the size of $n,m$ for (\ref{the_eq}) to have a unique solution). Finally we let 
\begin{equation}\label{neil}
\tilde f = \sum_{k=1}^n \tilde \beta_k w_k.
\end{equation}
Note that, for $n = m$ this is equivalent to (\ref{f}), and thus we have simply extended the framework discussed 
in Section \ref{early_work}. However, for $m>n$ this is no longer the case.  As we later establish, allowing $m$ to range independently of $n$ is the key to the advantage possessed by this framework.  

Before doing so, however, we first mention that the framework proposed above differs from that discussed previously in that it is inconsistent.  Unlike (\ref{f}), the samples 
$\langle s_{j},\tilde{f} \rangle$ do not coincide with those of the function $f$.  Yet, as we shall now see, by dropping the requirement of consistency, we obtain a reconstruction which circumvents the aforementioned issues associated with (\ref{f}).

\subsection{The Abstract Sampling Theorem} 
The task is now to analyze the model in (\ref{the_eq}) by both establishing existence of $\tilde{f}$ and providing error bounds for 
$\|f-\tilde f\|$.  We have
\begin{theorem}\label{abstract}
Let $\mathcal{H}$ be a separable Hilbert space and 
$\mathcal{S}, \mathcal{W} \subset \mathcal{H}$ be closed 
subspaces such that $\mathcal{W} \cap\mathcal{S}^{\perp} = \{0\}$. Suppose that $\{s_k\}_{k \in \mathbb{N}}$ and 
$\{w_k\}_{k \in \mathbb{N}}$ are Riesz bases for
$\mathcal{S}$ and $\mathcal{W}$ respectively with constants $A, B, C, D > 0$. 
Suppose that
\begin{equation}\label{fff}
f = \sum_{k\in\mathbb{N}} \beta_k w_k, \qquad \beta = \{\beta_1, \beta_2, \hdots, \} \in l^2(\mathbb{N}).  
\end{equation}
Let $n \in \mathbb{N}$. Then there is an $M \in \mathbb{N}$ (in particular 
$M = \min\{k: 0 \notin \sigma(P_nU^*P_kUP_n\lvert_{P_n\mathcal{H}}) \}$) such that, for all $m \geq M$, the solution  
$\{\tilde \beta_1,\hdots, \tilde \beta_n\}$ to (\ref{the_eq}) is unique. Also, if $\tilde f$ is as in 
(\ref{neil}), then 
\begin{equation}\label{err_b}
\|f - \tilde f\|_{\mathcal{H}} \leq \sqrt{B}(1+K_{n,m})\|P_n^{\perp}\beta\|_{l^2(\mathbb{N})}, 
\end{equation}
where 
\begin{equation}\label{Knm}
K_{n,m} = \left\|(P_nU^*P_mUP_n\lvert_{P_n\mathcal{H}})^{-1}P_nU^*P_mUP_n^{\perp}\right\|.
\end{equation}
\end{theorem}

The theorem has an immediate corollary that is useful for estimating the error.  We have
\begin{corollary}\label{abs_cor}
With the same assumptions as in Theorem \ref{abstract} and fixed $n \in \mathbb{N}$, 
\begin{equation}\label{con_m}
\left\|(P_nU^*P_mUP_n\lvert_{P_n\mathcal{H}})^{-1}\right\| \longrightarrow \left\|(P_nU^*UP_n\lvert_{P_n\mathcal{H}})^{-1}\right\| \leq \left\|(U^*U)^{-1}\right\| \leq \frac{1}{AC}, \quad m \rightarrow \infty.
\end{equation}
In addition, if  $U$ is an isometry (in particular, when 
$\{w_k\}_{k\in \mathbb{N}}, \{s_k\}_{k\in \mathbb{N}}$ are orthonormal) then it follows that
$$
K_{n,m} \longrightarrow 0, \qquad m \rightarrow \infty.
$$ 
\end{corollary}

\begin{proof}[Proof of Theorem \ref{abstract}]
Let $U$ be as in as in (\ref{inf_U}). Then (\ref{fff}) yields the following infinite system of equations:
\begin{equation}\label{recall44}
\left(
\begin{matrix}
 \langle s_1,f\rangle \\
 \langle s_2,f\rangle  \\
 \langle s_3,f\rangle  \\
\vdots       \\
\end{matrix}
\right)
= 
\left(
\begin{matrix}
 u_{11}         & u_{12}    & u_{13}  & \hdots \\
 u_{21}         & u_{22}    & u_{23}  & \hdots \\
 u_{31}         & u_{32}    & u_{33}  & \hdots \\
\vdots          & \vdots    & \vdots  & \ddots\\
\end{matrix}
\right) 
\left(
\begin{matrix}
 \beta_1 \\
 \beta_2 \\
 \beta_3 \\
\vdots       \\
\end{matrix}
\right).
\end{equation}
Note that $U$ must be a bounded operator. Indeed, let $S$ and $W$ be as in (\ref{SandW}).  Since 
 $$
\langle S^*W e_j,e_i \rangle = \langle s_i, w_j\rangle, \qquad i,j \in \mathbb{N},
 $$
 it follows that $U = S^*W$. However, from (\ref{riesz}) we find that both $W$ and $S$ are bounded as mappings from 
 $l^2(\mathbb{N})$ onto $\mathcal{W}$ and $\mathcal{S}$ respectively, with $\|W\| \leq \sqrt{B}$, 
 $\|S\| \leq \sqrt{D}$, thus yielding our claim.
Note also that, by the assumption that $\mathcal{W} \cap\mathcal{S}^{\perp} = \{0\}$, (\ref{recall44}) has a 
unique solution. Indeed,  since $\mathcal{W} \cap \mathcal{S}^{\perp} = \{0\}$ it follows that 
$  \inf_{\|x\|=1}\|S^*Wx\| \neq 0$, so $U$ must be injective.

 Now
 let $\eta_f = \{\langle s_1,f \rangle,  \langle s_1,f\rangle, \hdots\}$, then (\ref{recall44}) gives us that 
\begin{equation}\label{purple4}
P_nU^*P_m \eta_f = P_nU^*P_mU\left(P_n + P_n^{\perp}\right)\beta.
\end{equation}
Suppose for a moment that we can show that there exists an $M > 0$ such that 
$ P_nU^*P_mUP_n\lvert_{P_n\mathcal{H}}$ is invertible for all $m \geq M$.
Hence, we may appeal to (\ref{purple4}), whence
\begin{equation}\label{ups}
\begin{split}
(P_nU^*P_mUP_n\lvert_{P_n\mathcal{H}})^{-1}P_nU^*P_m \eta_f
 = P_n\beta + (P_nU^*P_mUP_n\lvert_{P_n\mathcal{H}})^{-1}P_nU^*P_mUP_n^{\perp}\beta,
\end{split}
\end{equation}
and therefore, by (\ref{purple4}) and (\ref{riesz}), 
\begin{equation*}
\begin{split}
\left\|f -  \sum_{k=1}^n \tilde \beta_k w_k\right\|_{\mathcal{H}}
&\leq
\sqrt{B} \left\|(P_nU^*P_mUP_n\lvert_{P_n\mathcal{H}})^{-1}
P_nU^*P_m \eta_f - \beta\right\|_{l^2(\mathbb{N})} \\ 
&= \sqrt{B}\left\|(P_n^{\perp}-  
(P_nU^*P_mUP_n\lvert_{P_n\mathcal{H}})^{-1}P_nU^*P_mUP_n^{\perp})
\beta\right\|_{l^2(\mathbb{N})} \\
&\leq \sqrt{B} \left(1 + K_{n,m} \right)
\left\|P_n^{\perp}\beta\right\|_{l^2(\mathbb{N})},\\
\end{split}
\end{equation*}
where 
$$
K_{n,m} = \left\|(P_nU^*P_mUP_n\lvert_{P_n\mathcal{H}})^{-1}
P_nU^*P_mUP_n^{\perp}\right\|.
$$
Thus, (\ref{err_b}) is established, provided we can show the following claim:

{\bf Claim:} There exists an $M > 0$ such that  
$P_nU^*P_mUP_n\lvert_{P_n\mathcal{H}}$ is invertible for all $m \geq M$. 
Moreover, 
$$
\left\|(P_nU^*P_mUP_n\lvert_{P_n\mathcal{H}})^{-1}\right\| \longrightarrow \left\|(P_nU^*UP_n\lvert_{P_n\mathcal{H}})^{-1}\right\| \leq \left\|(U^*U)^{-1}\right\|, \qquad m \rightarrow \infty.
$$

To prove the claim, we first need to show that $P_nU^*UP_n\lvert_{P_nl^2(\mathbb{N})}$ is invertible 
 for all $n \in \mathbb{N}$.   To see this, let $\Theta: \mathcal{B}(l^2(\mathbb{N})) \rightarrow \mathbb{C}$ denote 
 the numerical range.  Note that $U^*U$ is self-adjoint and invertible. The latter implies that there is a 
 neighborhood 
 $\omega$ around zero such that $\sigma(U^*U) \cap \omega = \emptyset $ and the former implies that the numerical range $\Theta(U^*U) \cap
\omega = \emptyset.$ 
Now the spectrum $\sigma(P_nU^*UP_n\lvert_{P_nl^2(\mathbb{N})}) 
 \subset \Theta(P_nU^*UP_n\lvert_{P_nl^2(\mathbb{N})}) \subset \Theta(U^*U).$ Thus, 
 $$
 \sigma(P_nU^*UP_n\lvert_{P_nl^2(\mathbb{N})})  \cap \omega = \emptyset, \qquad \forall 
 \, n \in \mathbb{N},
 $$ thus, 
 $P_nU^*UP_n\lvert_{P_nl^2(\mathbb{N})}$ 
 is always invertible. Now,
make the following two observations  
\begin{equation}\label{obs}
\begin{split}
P_nU^*P_m UP_n 
&= \sum_{j=1}^m (P_n \xi_j) \otimes (P_n \bar\xi_j), 
\qquad \xi_j = U^*e_j,\\
P_nU^*UP_n 
&= \sum_{j=1}^{\infty} (P_n \xi_j) \otimes (P_n \bar\xi_j), 
\end{split}
\end{equation}
where the last series converges at least strongly (it converges in norm, but that is a part of the proof).
The first is obvious.  The second observation follows from the fact that 
$P_m U \rightarrow U$ strongly as $m \rightarrow \infty.$
Note that
$$
\|P_n\xi_j\|^2 = \langle P_n\xi_j, P_n\xi_j \rangle = \langle U P_n U^*e_j, e_j, \rangle.
$$
However, $U^* P_nU$ must be trace class since $\mathrm{ran}(P_n)$ is finite-dimensional. Thus, by (\ref{obs}) 
we find that
\begin{equation*}
\begin{split}
 \|P_nU^*P_m UP_n - P_nU^*UP_n\| &\leq\sum_{j=m+1}^\infty \left\|(P_n \xi_j) \otimes (P_n \bar\xi_j)\right\|  \\
 &\leq  \sum_{j=m+1}^\infty \langle  U P_n U^*e_j, e_j, \rangle \longrightarrow 0, 
 \qquad m \rightarrow \infty.
 \end{split}
 \end{equation*}
 Hence,  the claim follows (the fact that  $\left\|(P_nU^*UP_n\lvert_{P_n\mathcal{H}})^{-1}\right\| \leq \left\|(U^*U)^{-1}\right\|$ is clear from the observation that $U^*U$ is self-adjoint) and we are done.  
 \end{proof}
 
 \begin{proof}[Proof of Corollary \ref{abs_cor}]
Note that the claim in the proof of Theorem \ref{abstract} yields the first part of (\ref{con_m}), and the second part follows from the fact that  $U = S^*W$ (where $S,W$ are also defined in the proof of Theorem \ref{abstract}) and 
 (\ref{riesz}). Thus, 
 we are now left with the task of showing that $K_{n,m} \rightarrow 0$ as $m \rightarrow \infty$ when $U$ is 
 an isometry. Note that the assertion will follow, by (\ref{Knm}), if we can show that 
$$
\left\|P_nU^*P_mUP_n^{\perp}\right\| \longrightarrow 0, \qquad m 
\longrightarrow \infty.
$$
However, this is straightforward, since a simple calculation yields
\begin{equation}\label{crep}
\left\|P_nU^*P_mUP_n^{\perp}\right\| \leq  \|U\|(\|P_nU^*P_m UP_n - P_nU^*UP_n\|)^{1/2}.
\end{equation}
\end{proof}

\begin{remark}
Note that the trained eye of an operator theorist will immediately spot that the claim in the proof of Theorem \ref{abstract} and Corollary \ref{abs_cor}  follows (with an easy reference to known convergence properties of finite rank operators in the 
strong operator topology) without the computations done in our exposition. However, we feel that the exposition 
illustrates ways of estimating bounds for 
$$
\left\|\left(P_nU^*P_mUP_n\lvert_{P_n\mathcal{H}}\right)^{-1}\right\|, \qquad \left\|P_nU^*P_mUP_n^{\perp}\right\|,
$$
which are crucial in order to obtain a bound for $K_{n,m}$.  This is demonstrated in Section \ref{Normbounds}.
\end{remark}

\begin{remark}
Note that $S^*W$ (and hence also $U$) is invertible if and only if $\mathcal{H} = \mathcal{W} \oplus \mathcal{S}^{\perp}$,
which is equivalent to  $\mathcal{W} \cap \mathcal{S}^{\perp} = \{0\}$ and $\mathcal{W}^{\perp} \cap \mathcal{S} = \{0\}$. This requirement is quite strong as we may very well have that $\mathcal{W} \neq \mathcal{H}$ and 
$\mathcal{S} = \mathcal{H}$ (e.g. Example 3.2 when $\epsilon<1$). In this case we obviously have that 
 $\mathcal{W}^{\perp} \cap \mathcal{S} \neq\{0\}$. However, as we saw in Theorem \ref{abstract}, as long as we have $f \in \mathcal{W}$ we only need injectivity of $U$, which is guaranteed when $\mathcal{W} \cap \mathcal{S}^{\perp} = \{0\}$.   
\end{remark}

If one wants to write our framework in the language used in Section \ref{early_work}, it is easy to see that our reconstruction 
can be written as  
\begin{equation}\label{our_f}
\tilde f  = W_n(W_n^*S_mS_m^*W_n)^{-1}W_n^*S_mS_m^*f,
\end{equation}
where the operators $S_m : \mathbb{C}^m \rightarrow \mathcal{H}$ and 
$W_n : \mathbb{C}^n \rightarrow \mathcal{H}$ are defined as in (\ref{S_and_W}), and 
$S_m$ and $W_n$ corresponds to the spaces 
\begin{equation}\label{spaces}
\mathcal{S}_m = \mathrm{span}\{s_1,\hdots,s_m\}, \qquad \mathcal{W}_n = \mathrm{span}\{w_1,\hdots,w_n\},
\end{equation} 
where $\{w_k\}_{k\in\mathbb{N}}$ 
and $\{s_k\}_{k\in\mathbb{N}}$ are as in Theorem \ref{abstract}.
In particular, we get the following corollary:

\begin{corollary}
Let $\mathcal{H}$ be a separable Hilbert space and 
$\mathcal{S}, \mathcal{W} \subset \mathcal{H}$ be closed 
subspaces such that $\mathcal{W} \cap\mathcal{S}^{\perp} = \{0\}$. Suppose that $\{s_k\}_{k \in \mathbb{N}}$ and $\{w_k\}_{k \in \mathbb{N}}$ are Riesz bases for
$\mathcal{S}$ and $\mathcal{W}$ respectively. 
Then, for each $n \in \mathbb{N}$ there is an $M \in \mathbb{N}$ 
such that, for all $m \geq M$, the mapping $W_n^*S_mS_m^*W_n:\mathbb{C}^n \rightarrow \mathbb{C}^n$ is invertible 
(with $S_m$ and $W_n$ defined as above). Moreover, if $\tilde f$ is as in (\ref{our_f}), then
$$
\left\|P_{\mathcal{W}_n}^{\perp}f \right\|_{\mathcal{H}} \leq \|f-\tilde f\|_{\mathcal{H}} \leq (1+ K_{n,m})\left\|P_{\mathcal{W}_n}^{\perp}f \right\|_{\mathcal{H}},
$$
where $P_{\mathcal{W}_n}$ is the orthogonal projection onto $\mathcal{W}_n$, and
$$
K_{n,m} = \left\|W_n(W_n^*S_mS_m^*W_n)^{-1}W_n^*S_mS_m^*P_{\mathcal{W}_n}^{\perp}\right\|. 
$$
Moreover, when $\{ s_k \}$ and $\{ w_k \}$ are orthonormal bases, then, for fixed $n$, $C_{n,m} \rightarrow 0$ as $m\rightarrow \infty$. 
\end{corollary}

\begin{proof}
The fact that $W_n^*S_mS_m^*W_n:\mathbb{C}^n \rightarrow \mathbb{C}^n$ is invertible for large $m$ follows from the the observation that $\mathcal{W} \cap\mathcal{S}^{\perp} = \{0\}$ and the proof 
of Theorem \ref{abstract}, by noting that $S_m^*W_n = P_mUP_n$, where $U$ is as in Theorem \ref{abstract}.
Now observe that 
\begin{equation}\label{univ}
\begin{split}
W_n^*S_mS_m^*f &= W_n^*S_mS_m^*(P_{\mathcal{W}_n}f + P_{\mathcal{W}_n}^{\perp}f )\\
&=  W_n^*S_mS_m^* W_n(W^*_nW_n)^{-1}W_n^*f + W_n^*S_mS_m^* P_{\mathcal{W}_n}^{\perp}f. 
\end{split}
\end{equation}
Note also that $W^*_nW_n:\mathbb{C}^n \rightarrow \mathbb{C}^n$ is clearly invertible, since 
$\{w_k\}_{k =1}^n$ are linearly independent. Now (\ref{univ}) yields
$$
W_n(W_n^*S_mS_m^*W_n)^{-1}W_n^*S_mS_m^*f = P_{\mathcal{W}_n}f  + W_n(W_n^*S_mS_m^*W_n)^{-1}W_n^*S_mS_m^*P_{\mathcal{W}_n}^{\perp}f.
$$
Thus, 
$$
\|f-\tilde f\|_{\mathcal{H}} \leq \left\|P_{\mathcal{W}_n}^{\perp}  -  W_n(W_n^*S_mS_m^*W_n)^{-1}W_n^*S_mS_m^*P_{\mathcal{W}_n}^{\perp}\right\|_{\mathcal{H}} \left\|P_{\mathcal{W}_n}^{\perp}f \right\|_{\mathcal{H}},
$$
which gives the first part of the corollary. The second part follows from similar reasoning as in the proof of Corollary \ref{abs_cor}. 
\end{proof}

\begin{remark}
The framework explained in Section 
\ref{early_work} is equivalent to using the finite section method, and this may work for certain bases, however, not in general (as Example \ref{exp} shows).  Computing with infinite matrices can be a challenge since the qualities of 
any finite section may be very different from the original infinite matrix. The use of uneven sections (as we do in this paper) of infinite matrices seems to be the best way to combat these problems. The reader may consult \cite{strohmer,hansen2011, Lindner2008} for other examples of uneven section techniques. 
\end{remark}

When compared to the method of Eldar et al, the framework presented here has a number of important advantages:
\begin{itemize}
\item[(i)] It allows reconstructions in arbitrary bases and does not need extra assumptions as in (\ref{SW}).
\item[(ii)]  The conditions on $m$ (as a function of $n$) for $P_nU^*P_mUP_n\lvert_{P_n\mathcal{H}}$ to be invertible (such that we have a unique solution) can be numerically computed. Moreover, bounds on the 
constant $K_{n,m}$ can also be computed efficiently. This is the topic in Section \ref{Normbounds}.  
\item[(iii)] It is numerically stable: the matrix $A = P_n U^{*} P_m U P_n |_{P_n \mathcal{H}}$ has bounded inverse (Corollary \ref{abs_cor}) for all $n$ and $m$ sufficiently large.
\item[(iv)] The approximation $\tilde{f}$ is quasi-optimal (in $n$).  It converges at the same rate as the tail $\| P^{\perp}_{n} \beta \|_{l^2(\mathbb{N})}$, in contrast to (\ref{f}) which converges more slowly whenever the parameter $\frac{1}{\cos ( \theta_{\mathcal{W}_m \mathcal{S}_m})}$ grows with $n=m$.
\end{itemize}
As mentioned, this method is inconsistent.  However, since $\{ s_{j} \}$ is a Riesz basis, we deduce that
\begin{equation*}
\sum^{m}_{j=1} | \langle  s_{j}, f - \tilde{f}  \rangle |^2 \leq c \| f - \tilde{f} \|^2,
\end{equation*}
for some constant $c>0$.  Hence, the departure from consistency (i.e. the left-hand side) is bounded by a constant multiple of the approximation error, and thus can also be bounded by $\| P^{\perp}_{n} \beta \|_{l^2(\mathbb{N})}$.

\subsection{The Generalized (Nyquist-Shannon) Sampling Theorem}\label{Sampling}

In this section, we apply the abstract sampling theorem (Theorem \ref{abstract}) 
to the classical sampling problem of recovering a function from samples of its Fourier transform.  As we shall see, when considered in this way, the corresponding theorem, which we call the generalized (Nyquist--Shannon) Sampling Theorem, extends the classical Shannon theorem (which corresponds to a special case) by allow reconstructions in arbitrary bases.

\begin{proposition}\label{prop}
Let $\mathcal{F}$ denote the Fourier transform on $L^2(\mathbb{R}^d).$
Suppose that $\{\varphi_j\}_{j \in \mathbb{N}}$ is a Riesz basis with constants $A,B$ 
(as in (\ref{riesz})) for a subspace 
$\mathcal{W} \subset L^2(\mathbb{R}^d)$ such that there exists a $T > 0$ with
$\mathrm{supp}(\varphi_j) \subset [-T,T]^d$ for all $j \in
\mathbb{N}.$
For $\epsilon > 0$, let 
$\rho: \mathbb{N} \rightarrow (\epsilon \mathbb{Z})^d$ be a
bijection.  Define the infinite matrix
$$
U =
\left(
\begin{matrix}
 u_{11}         & u_{12}    & u_{13}  & \hdots \\
 u_{21}         & u_{22}    & u_{23}  & \hdots \\
 u_{31}         & u_{32}    & u_{33}  & \hdots \\
\vdots          & \vdots    & \vdots  & \ddots\\
\end{matrix}
\right), \qquad u_{ij} = (\mathcal{F}\varphi_j)(\rho(i)).
$$
Then, for $\epsilon \leq \frac{1}{2T}$, we have that $U: l^2(\mathbb{N}) \rightarrow l^2(\mathbb{N})$ is 
bounded and invertible on its range with $\|U\| \leq \sqrt{\epsilon^{-d}B}$ and 
$\|(U^*U)^{-1}\| \leq \epsilon^dA^{-1}$ . Moreover, if $\{\varphi_j\}_{j \in \mathbb{N}}$ 
is an orthonormal set, then $\epsilon^{d/2}U$ is an isometry.
\end{proposition}

\begin{theorem}(The Generalized Sampling Theorem)\label{Diff_Shan}
With the same setup as in Proposition \ref{prop}, set
$$
f = \mathcal{F}g, \quad g = \sum_{j=1}^{\infty} \beta_j \varphi_j  \in L^2(\mathbb{R}^d),
$$ and let $P_n$ denote 
the projection onto $\mathrm{span}\{e_1,\ldots,e_n\}$.
Then, for every $n \in \mathbb{N}$ there is an $M \in \mathbb{N}$ such that, for all $m \geq M$, the solution to
$$
A\left(
\begin{matrix}
 \tilde \beta_1 \\
 \tilde \beta_2 \\
\vdots       \\
\tilde \beta_n
\end{matrix}
\right)
=
P_nU^*P_m
\left(
\begin{matrix} 
 f(\rho(1)) \\
f(\rho(2))  \\
 \vdots       \\
f(\rho(m))  \\
\end{matrix}
\right), \qquad  A = P_nU^*P_mUP_n\lvert_{P_n\mathcal{H}},
$$
is unique. Also, if 
$$
\tilde g =   \sum_{j=1}^{n}\tilde \beta_j \varphi_j,  \qquad 
\tilde f =   \sum_{j=1}^{n}\tilde \beta_j \mathcal{F}\varphi_j,
$$
then 
\begin{equation}\label{first}
\|g - \tilde g \|_{L^2(\mathbb{R}^d)} \leq \sqrt{B}(1+K_{n,m})\|P_n^{\perp}\beta\|_{l^2(\mathbb{N})}, 
\quad \beta = \{\beta_1, \beta_2, \hdots\},
\end{equation}
and
\begin{equation}\label{second}
\|f - \tilde f\|_{L^{\infty}(\mathbb{R}^d)} \leq  (2T)^{d/2}  \sqrt{B}(1+K_{n,m})\|P_n^{\perp}\beta\|_{l^2(\mathbb{N})},
\end{equation}
where $K_{n,m}$ is given by (\ref{Knm}) and satisfies (\ref{con_m}). Moreover, when 
$\{\varphi_j\}_{j \in \mathbb{N}}$ is an orthonormal set, we have
$$
K_{n,m} \longrightarrow 0, \qquad m \rightarrow \infty,
$$
for fixed $n$.
\end{theorem}

\begin{proof}[Proof of Proposition \ref{prop}]
Note that 
\begin{equation*}
\begin{split}
u_{ij} = \int_{\mathbb{R}^d} \varphi_j(x) e^{-2\pi \mathrm{i} \rho(i) \cdot x}\, dx 
= \int_{[-T,T]^d} \varphi_j\left(x\right) e^{-2\pi \mathrm{i} \rho(i) \cdot x}\, dx.
\end{split}
\end{equation*}
Since $\rho: \mathbb{N} \rightarrow (\epsilon \mathbb{Z})^N$ is a bijection, it follows
that the functions $\{x \mapsto \epsilon^{d/2}e^{-2\pi \mathrm{i} \rho(i) \cdot x}\}_{i\in\mathbb{N}}$ form
an orthonormal basis for 
$L^2([-(2\epsilon)^{-1},(2\epsilon)^{-1}]^d) \supset L^2([-T,T]^d)$.
Let 
$$
\langle \cdot,\cdot \rangle = \overline{\langle \cdot,\cdot \rangle}_{L^2([-(2\epsilon)^{-1},(2\epsilon)^{-1}]^d)},
$$
denote a new inner product on $L^2([-(2\epsilon)^{-1},(2\epsilon)^{-1}]^d)$. 
Thus, we are now in the setting of Theorem \ref{abstract} and Corollary \ref{abs_cor} 
with $C = D =\epsilon^d$. It follows by Theorem 
\ref{abstract} and Corollary \ref{abs_cor} that $U$ is bounded and invertible on its range with $\|U\| \leq \sqrt{\epsilon^{-d}B}$ and $\|(U^*U)^{-1}\| \leq \epsilon^dA^{-1}$. Also, $\epsilon^{d/2}U$ is an isometry whenever $A = B = 1$, in particular when $\{\varphi_k\}_{k\in\mathbb{N}}$ is an orthonormal set.
\end{proof}

\begin{proof}[Proof of Theorem \ref{Diff_Shan}]
Note that (\ref{first}) now automatically follows from Theorem \ref{abstract}. To get (\ref{second}) we simply 
observe that, by the definition of the Fourier transform and using the Cauchy--Schwarz inequality,
\begin{equation*}
\begin{split}
&\sup_{x \in \mathbb{R}^d}\left|f(x) - \sum_{j=1}^n \tilde \beta_j \mathcal{F}\varphi_j(x)\right| \leq \int_{[-T,T]^d} \left|g(y) -  \sum_{j=1}^n \tilde \beta_j \varphi_j(y) \right| \, dy \\
&\leq 
(2T)^{d/2}\left \|g -  \sum_{j=1}^n \tilde \beta_j \varphi_j\right\|_{L^2(\mathbb{R}^d)} \leq
 (2T)^{d/2}\sqrt{B}(1+K_{n,m})\|P_n^{\perp}\beta\|_{l^2(\mathbb{N})},
\end{split}
\end{equation*}
where the last inequality follows from the already established (\ref{first}).  Hence we are done with the first part of 
the theorem. To see that $K_{n,m} \rightarrow 0$ as $m \rightarrow \infty$ when $\{\varphi_j\}_{j \in \mathbb{N}}$ 
is an orthonormal set, we observe that orthonormality yields $A=B=1$ and hence (since we already have established the values of $C$ and $D$) $\epsilon^{d/2}U$ must be 
an isometry. The convergence to zero now follows from Theorem \ref{abstract}. 
\end{proof}

Note that the bijection $\rho: \mathbb{N} \rightarrow (\epsilon \mathbb{Z})^d$ is only important when 
$d > 1$ to obtain an operator $U:l^2(\mathbb{N}) \rightarrow l^2(\mathbb{N})$.  However, when 
$d = 1,$ there is nothing preventing us from avoiding $\rho$ and forming an operator   
$U:l^2(\mathbb{N}) \rightarrow l^2(\mathbb{Z})$ instead. The idea follows below.
Let $\mathcal{F}$ denote the Fourier transform on $L^2(\mathbb{R}),$ and let $f = \mathcal{F}g$ for some 
$g \in L^2(\mathbb{R})$. Suppose that 
$\{\varphi_j\}_{j \in \mathbb{N}}$ is a Riesz basis for a closed subspace in $L^2(\mathbb{R})$ with constants 
$A,B>0$,
such that there is a $T > 0$ with
$\mathrm{supp}(\varphi_j) \subset [-T,T]$ for all $j \in \mathbb{N}.$ For $\epsilon > 0$, let 
\begin{equation}\label{hatU}
\widehat U =
\left(
\begin{matrix}
\vdots            & \vdots     & \vdots   &  \iddots\\
 u_{-1,1}         & u_{-1,2}    & u_{-1,3}  & \hdots \\
 u_{0,1}         & u_{0,2}    & u_{0,3}  & \hdots \\
 u_{1,1}         & u_{1,2}    & u_{1,3}  & \hdots \\
\vdots          & \vdots    & \vdots  & \ddots\\
\end{matrix}
\right), \qquad u_{i,j} = (\mathcal{F}\varphi_j)(i\epsilon).
\end{equation}
Thus, as argued in the proof of Theorem \ref{Diff_Shan}, $\widehat U \in \mathcal{B}(l^2(\mathbb{N}), l^2(\mathbb{Z})),$ 
provided $\epsilon \leq \frac{1}{2T}$. 
Next, let $P_n \in \mathcal{B}(l^2(\mathbb{N}))$ and, for odd $m$, $\tilde{P}_{m} \in \mathcal{B}(l^2(\mathbb{Z}))$ be the projections onto 
$$
\mathrm{span}\{e_1,\hdots, e_n\}, \qquad \mathrm{span}\{e_{-\frac{m-1}{2}},\hdots, e_{\frac{m-1}{2}}\}
$$ respectively.  Define $\{\tilde \beta_1, \hdots, \tilde \beta_n\}$ by 
(this is understood to be for sufficiently large $m$)
\begin{equation}\label{black}
\widehat A\left(
\begin{matrix}
 \tilde \beta_1 \\
 \tilde \beta_2 \\
 \tilde \beta_3 \\
\vdots       \\
\tilde \beta_n\\
\end{matrix}
\right)
=
P_n\widehat U^*P_m
\left(
\begin{matrix}
f(-\frac{m-1}{2}) \\
\vdots  \\
f(0)  \\
\vdots       \\
 f(\frac{m-1}{2})  \\
\end{matrix}
\right), 
\, \, \widehat A = P_n\widehat U^*P_m\widehat UP_n\lvert_{P_n\mathcal{H}}.
\end{equation}
By exactly the same arguments as in the proof 
of Theorem \ref{Diff_Shan}, it follows
that, if
$
g  = \sum_{j=1}^{\infty} \beta_j \varphi_j$, $\tilde g =   \sum_{j=1}^{n}\tilde \beta_j \varphi_j$, $f = \mathcal{F}g$ 
and $\tilde f =  \sum_{j=1}^{n}\tilde \beta_j\mathcal{F} \varphi_j$, then
\begin{equation}\label{black2}
\begin{split}
\|g - \tilde g \|_{L^2(\mathbb{R})}
 &\leq \sqrt{B}(1+K_{n,m})\|P_n^{\perp}\beta\|_{l^2(\mathbb{N})}, 
\,\, \beta = \{\beta_1, \beta_2, \hdots\},\\
\|f-\tilde f\|_{L^{\infty}(\mathbb{R})}
&\leq  \sqrt{2T}\sqrt{B}(1+K_{n,m})\|P_n^{\perp}\beta\|_{l^2(\mathbb{N})},
\end{split}
\end{equation}
where  $K_{n,m}$ is as in (\ref{Knm}).

\begin{remark}
Note that (as the proof of the next corollary will show) the classical NS-Sampling 
Theorem is just a special case of Theorem \ref{Diff_Shan}.
\end{remark}

\begin{corollary}
Suppose that $f = \mathcal{F}g$ and $\mathrm{supp}(g) \subset [-T,T]$. 
Then, for $ 0 < \epsilon \leq \frac{1}{2T}$ we have that
$$
g(\cdot) = \epsilon \sum_{k=-\infty}^{\infty}f(k\epsilon)e^{2\pi \mathrm{i} \epsilon k \cdot }
\qquad L²\text{ convergence}.
$$ 
$$
f(t) = \sum_{k=-\infty}^{\infty}f(k\epsilon)\mathrm{sinc}\left(\frac{t
    +k\epsilon}{\epsilon}\right) \qquad L²\,\text{and unif. convergence.} 
$$
\end{corollary}

\begin{proof}
Define the basis $\{\varphi_j\}_{j \in \mathbb{N}}$ for 
$L^2([-(2\epsilon)^{-1},(2\epsilon)^{-1}])$ by 
\begin{equation*}
\begin{split}
\varphi_1(x) = 
\sqrt{\epsilon}&\chi_{[-\frac{1}{2\epsilon},
  \frac{1}{2\epsilon}]}(x), \quad 
\varphi_2(x) = 
\sqrt{\epsilon}e^{2\pi \mathrm{i} \epsilon x}\chi_{[-\frac{1}{2\epsilon},
  \frac{1}{2\epsilon}]}(x), \\
\varphi_3(x) &= 
\sqrt{\epsilon}e^{2\pi \mathrm{i} \epsilon (-1)x} \chi_{[-\frac{1}{2\epsilon},
  \frac{1}{2\epsilon}]}(x), \\
\varphi_4(x) &= 
\sqrt{\epsilon}e^{2\pi \mathrm{i} \epsilon 2x}\chi_{[-\frac{1}{2\epsilon},
  \frac{1}{2\epsilon}]}(x), \\ 
 \varphi_5(x) &= 
\sqrt{\epsilon}e^{2\pi \mathrm{i} \epsilon (-2)x} \chi_{[-\frac{1}{2\epsilon},
  \frac{1}{2\epsilon}]}(x), \\
\varphi_{6}(x) &= 
\sqrt{\epsilon}e^{2\pi \mathrm{i} \epsilon 3 x}\chi_{[-\frac{1}{2\epsilon},
  \frac{1}{2\epsilon}]}(x)\quad \text{etc.} 
\end{split}
\end{equation*}
Letting $\widehat U = \{u_{k,l}\}_{k \in \mathbb{Z}, l \in \mathbb{N}},$ where 
$
u_{k,l} = (\mathcal{F}\varphi_l)(k\epsilon), 
$ an easy computation shows that  
$$
\widehat U =
\left(
\begin{matrix}
\vdots     & \vdots      & \vdots     & \vdots   & \vdots& \iddots\\
 0         & 0    & 0  & 0 & \frac{1}{\sqrt{\epsilon}} & \hdots \\
 0         & 0    & \frac{1}{\sqrt{\epsilon}}  & 0 & 0 & \hdots \\
 \frac{1}{\sqrt{\epsilon}}         & 0    & 0  & 0 & 0 & \hdots \\
 0         & \frac{1}{\sqrt{\epsilon}}    & 0  & 0 & 0 & \hdots \\
 0         & 0    & 0  & \frac{1}{\sqrt{\epsilon}} & 0 & \hdots \\
\vdots     & \vdots      & \vdots    &  \vdots  & \vdots &\ddots\\
\end{matrix}
\right).
$$
By choosing $m = n$ in (\ref{black}), we find that  $\tilde \beta_1 = \sqrt{\epsilon}f(0),$ $\tilde \beta_2 =  \sqrt{\epsilon}f(\epsilon),$ 
$\tilde \beta_3 =  \sqrt{\epsilon}f(-\epsilon),$ etc and that $K_{n,m} = 0$ in (\ref{black2}).
The corollary then follows from (\ref{black2}). 
\end{proof}

\begin{remark}
Returning to the general case, recall the definition of $\Omega_{N,\epsilon} $ from (\ref{one}), the mappings 
$\Lambda_{N,\epsilon,1}$, $\Lambda_{N,\epsilon,2}$ from (\ref{two}) and 
$\Theta$ from (\ref{three}). 
Define $\Xi_{N,\epsilon,1}: \Omega_{N,\epsilon} \rightarrow L^2(\mathbb{R})$ and 
$\Xi_{N,\epsilon,2}: \Omega_{N,\epsilon} \rightarrow L^2(\mathbb{R})$ by
$$
\Xi_{N,\epsilon,1}(f) = \sum_{j=1}^n \tilde \beta_j \mathcal{F}\varphi_j(\cdot), \qquad \Xi_{N,\epsilon,2}(f) = \sum_{j=1}^n \tilde \beta_j \varphi_j(\cdot),
$$
where $\tilde \beta = \{\tilde \beta_1, \hdots, \tilde \beta_n\}$ is the solution to (\ref{black}) with $N = m.$
Then, for $n > M$ (recall $M$ from the definition of $\Theta$ (\ref{three})), and
$$ 
m =  m(\gamma) = 
\min \{k \in \mathbb{N}: \|(P_n\widehat U^*P_k\widehat UP_n\lvert_{P_n\mathcal{H}})^{-1}\| \leq \epsilon \gamma\},\quad \gamma>1,
$$
it follows that 
\begin{equation*}
\begin{split}
\|\Xi_{N,\epsilon,1}(f) - f\|_{L^{\infty}(\mathbb{R})} &= 0 <
\|\Lambda_{N,\epsilon,1}(f) - f\|_{L^{\infty}(\mathbb{R})} 
\forall f, f = \mathcal{F}g, g \in \Theta,\\
\|\Xi_{N,\epsilon,2}(f) - g\|_{L^{2}(\mathbb{R})} &= 0 < 
\|\Lambda_{N,\epsilon,2}(f) - g\|_{L^{2}(\mathbb{R})} 
 \forall f, f = \mathcal{F}g, g \in \Theta.
\end{split}
\end{equation*}
Hence, under the aforementioned assumptions on $m$ and $n$, both $f$ and $g$ are recovered exactly by this method, provided $g \in \Theta$. Moreover, the reconstruction is done in a 
stable manner, where the stability depends on the parameter $\gamma$.
\end{remark}
To complete this section, let us sum up several of the key features of Theorem \ref{Diff_Shan}.  First, whenever $m$ is sufficiently large, the error incurred by $\tilde{g}$ is directly related to the properties of $g$ with respect to the reconstruction basis.  In particular, as previously noted, $g$ is reconstructed exactly under certain conditions.  Second, for fixed $n$, by increasing $m$ we can get arbitrarily close to the best approximation to $g$ in the reconstruction basis whenever the reconstruction vectors are orthonormal.  Thus, provided an appropriate basis is known, this procedure allows for near-optimal recovery.  The main question that remains, however, is how to guarantee that the conditions of Theorem \ref{Diff_Shan} are satisfied.  This is the topic of the next section.

\section{Norm Bounds}\label{Normbounds}

\subsection{Determining $m$}\label{m}
Recall that the constant $K_{n,m}$ in the error bound in Theorem \ref{abstract} (recall also $U$ from the same theorem) is given by
$$
K_{n,m} = \left\|(P_nU^*P_mUP_n\lvert_{P_n\mathcal{H}})^{-1}P_nU^*P_mUP_n^{\perp}\right\|.
$$
It is therefore of utmost importance to estimate $K_{n,m}$. This can be done numerically. Note that we already 
have established bounds on $\|U\|$ depending on the Riesz constants in (\ref{riesz}) and since we obviously have 
that $K_{n,m} \leq \|(P_nU^*P_mUP_n\lvert_{P_n\mathcal{H}})^{-1}\|\|U\|^2,$ we only require an estimate for the quantity 
$\|(P_nU^*P_mUP_n\lvert_{P_n\mathcal{H}})^{-1}\|$.

Recall also from Theorem \ref{abstract} that, if $U$ is an isometry up to a constant, then $K_{n,m} \rightarrow 0$ 
as $m \rightarrow \infty.$ In the rest of this section we will assume that $U$ has this quality.
In this case we are interested in the following problem:
given $n \in \mathbb{N}, \theta \in \mathbb{R}_+$, what is the smallest $m \in \mathbb{N}$ such that 
$K_{n,m} \leq \theta$? More formally, we wish to estimate the function 
$\Phi: \mathcal{U}(l^2(\mathbb{N})) \times \mathbb{N} \times \mathbb{R}_+ \rightarrow \mathbb{N},$ 
\begin{equation}\label{thefunction}
\begin{split}
\Phi(U,n,\theta) = \min\left\{m \in \mathbb{N}:
\left\|(P_nU^*P_mUP_n\lvert_{P_n\mathcal{H}})^{-1}P_nU^*P_mUP_n^{\perp}\right\| \leq \theta\right\},
\end{split}
\end{equation} 
where 
$$
\mathcal{U}(l^2(\mathbb{N})) = \left\{U \in \mathcal{B}(l^2(\mathbb{N})): U^*U = cI, c \in \mathbb{R}_+\right\}.
$$
Note that $\Phi$ is well defined for all $\theta \in \mathbb{R}_+$, since we have established that 
$K_{n,m} \rightarrow 0$ as $m \rightarrow \infty.$

\subsection{Computing Upper and Lower Bounds on $K_{n,m}$}
\begin{figure}
\centering
\includegraphics[height=29mm]{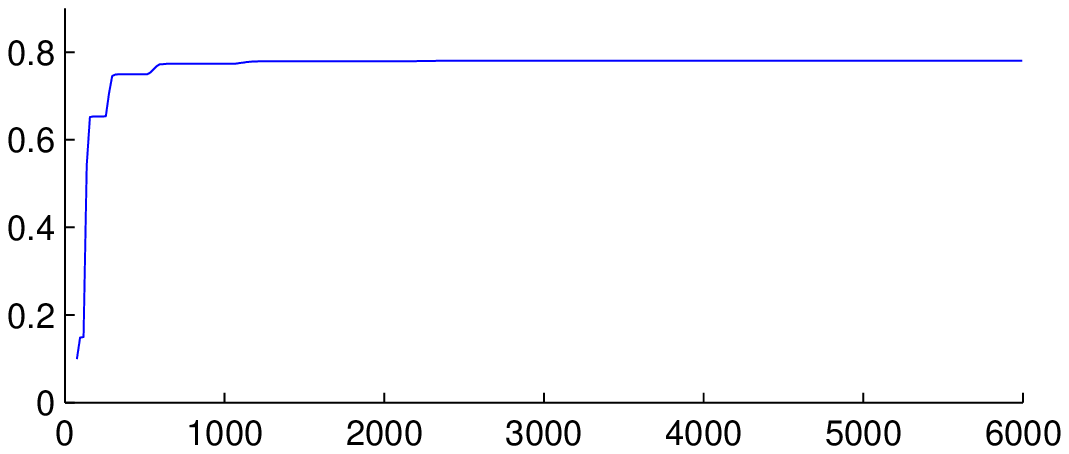}
\includegraphics[height=29mm]{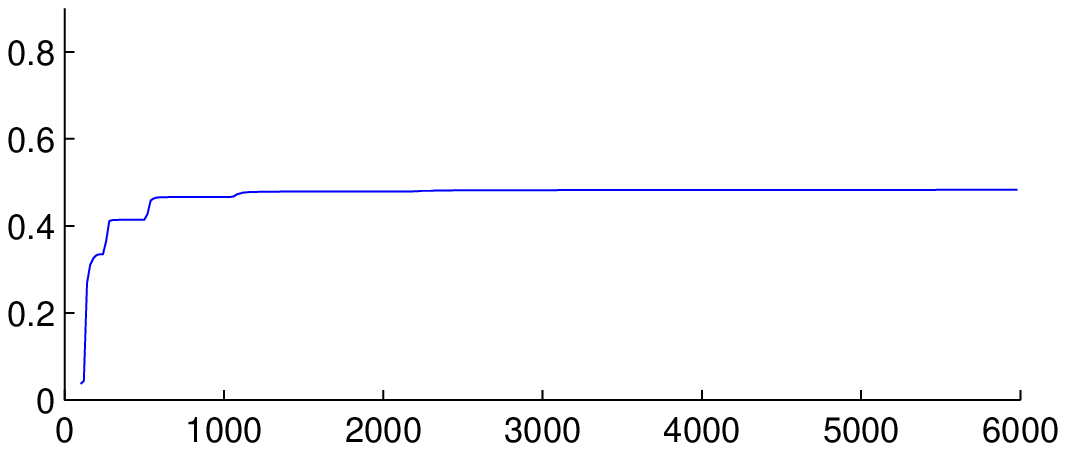}
\caption{The figure shows $K_{n,m,M}$ for $n = 75$, $m=350$ and $M = n+1, \hdots, 6000$ (left) and 
$K_{n,m,M}$ for $n = 100$, $m=400$ and $M = n+1, \hdots, 6000$ (right) for the Haar wavelets on $[0,1]$.  }
\label{fred2}
\end{figure}

The fact that $UP_n^{\perp}$ has infinite rank makes the computation of $K_{n,m}$ 
a challenge.  However, we may compute approximations from above and below.  For $M \in \mathbb{N}$, define
$$
K_{n,m,M} =  \left\|(P_nU^*P_mUP_n \lvert_{P_n\mathcal{H}})^{-1}P_nU^*P_mUP_n^{\perp}P_M\right\|,
$$
$$
 \widetilde K_{n,m} =  \left\|(P_nU^*P_mUP_n\lvert_{P_n\mathcal{H}})^{-1}P_nU^*P_m\right\|.
$$ 
Then,  for $L \geq M$
\begin{equation*}
\begin{split}
K_{n,m,M}
 &= \sup_{\xi \in P_M\mathcal{H}, \|\xi\| = 1}\left\|(P_nU^*P_mUP_n\lvert_{P_n\mathcal{H}})^{-1}P_nU^*P_mUP_n^{\perp}P_M\xi\right\| \\
&\leq \sup_{\xi \in P_L\mathcal{H}, \|\xi\| = 1}\left\|(P_nU^*P_mUP_n\lvert_{P_n\mathcal{H}})^{-1}P_nU^*P_mUP_n^{\perp}P_L\xi\right\|\\
&\leq  \sup_{\xi \in \mathcal{H}, \|\xi\| = 1}\left\|(P_nU^*P_mUP_n\lvert_{P_n\mathcal{H}})^{-1}P_nU^*P_mUP_n^{\perp}\xi\right\|
= K_{n,m}.
\end{split}
\end{equation*}
Clearly, $ K_{n,m} \leq   \|U\| \widetilde K_{n,m}$ and, since $P_M\xi \rightarrow \xi$ as $M \rightarrow \infty$ for 
all $\xi \in \mathcal{H}$, and by the reasoning above, it follows that 
$$
K_{n,m,M} \leq K_{n,m} \leq \|U\| \widetilde K_{n,m}, \qquad  K_{n,m,M} \nearrow K_{n,m}, \quad M \rightarrow \infty. 
$$
Note that 
$$
(P_nU^*P_mUP_n\lvert_{P_n\mathcal{H}})^{-1}P_nU^*P_mUP_n^{\perp}P_M: 
P_M \mathcal{H} \rightarrow P_n\mathcal{H}
$$
has finite rank. Therefore we may easily compute $K_{n,m,M}$. In Figure \ref{fred2} we have computed 
$K_{n,m,M}$ for different values of $n,m,M$. Note the rapid convergence in both examples.

\subsection{Wavelet bases}
Whilst in the general case $\Phi(U,n,\theta)$ must be computed numerically, in certain cases we are able to derive explicit analytical bounds for this quantity.  As an example, we now describe how to obtain bounds for bases consisting of compactly supported wavelets.  Wavelets and their various generalizations present an extremely efficient means in which to represent functions (or, more commonly, signals) \cite{Cand, Candes2002, Mallat1998}.  Given their long list of applications, the development of wavelet-based reconstruction methods is naturally a topic of utmost importance.

Let us review the basic wavelet approach on how to create orthonormal subsets 
$\{\varphi_k\}_{k\in \mathbb{N}} \subset L^2(\mathbb{R})$ with the property that 
$L^2([0,a]) \subset \mathrm{cl}(\mathrm{span}\{\varphi_k\}_{k\in\mathbb{N}})$ for some $a>0$.  Suppose that we are given a mother wavelet $\psi$ and a scaling function $\phi$ such that 
$\mathrm{supp}(\psi) = \mathrm{supp}(\phi) = [0,a]$ for some $a>0$.
The most obvious approach is to consider the following collection of functions:
$$
\Omega_a = \{\phi_k, \psi_{j,k}: j \in \mathbb{Z}_+, k \in \mathbb{Z}, \mathrm{supp}(\phi_k)^o \cap [0,a] \neq \emptyset, \, \mathrm{supp}(\psi_{j,k})^o \cap [0,a] \neq \emptyset \},
$$
where 
$$
\phi_k = \phi(\cdot - k), \qquad \psi_{j,k} = 2^{\frac{j}{2}}\psi(2^j \cdot - k).
$$
(The notation $K^o$ denotes the interior of a set $K \subset \mathbb{R}$.) Then we will have that 
$$
L^2([0,a]) \subset \mathrm{cl}(\mathrm{span}\{\varphi: \varphi \in \Omega_a\}) \subset L^2[-T,T],
$$
where $T > 0$ is such that $[-T,T]$ contains the support of all functions in $\Omega_a$.
However, the inclusions may be proper (but not always, as is the case with the Haar wavelet.)
It is easy to see that 
$$
 \psi_{j,k} \notin \Omega_a   \Longleftrightarrow \frac{a+k}{2^j} \leq 0, \quad a \leq \frac{k}{2^j},
$$
$$
\phi_k \notin \Omega_a   \Longleftrightarrow a+k \leq 0, \quad a \leq k.
$$
Hence we get that 
$$
\Omega_a = \{\phi_k: |k| = 0,\hdots, \lceil a\rceil-1\} \cup \{\psi_{j,k}: j \in \mathbb{Z}_+, k \in \mathbb{Z},   -\lceil a\rceil+1 \leq k \leq 2^j\lceil a\rceil-1\},
$$
and we will order $\Omega_a$ as follows:
\begin{equation}\label{med}
\{\phi, \phi_1, \hdots, \phi_{\lceil a\rceil-1}, \phi_{-1}, \hdots, \phi_{-\lceil a\rceil+1}, \psi_{0,0},  \psi_{0,1}, \hdots, \psi_{0,\lceil a\rceil-1}, \psi_{0,-1} ,\hdots,\psi_{0,-\lceil a\rceil+1}, \psi_{1,0}, \hdots \}.
\end{equation}
We will in this section be concerned with compactly supported wavelets and scaling functions satisfying
\begin{equation}\label{wavelet_C}
|\mathcal{F}\phi(w)| \leq \frac{C}{|w|^p}, \qquad |\mathcal{F}\psi(w)| \leq \frac{C}{|w|^p}, \qquad \omega \in \mathbb{R}\setminus \{0\},
\end{equation}
for some 
$$
C > 0, \qquad p \in \mathbb{N}.
$$
Before we state and prove bounds on $\Phi(U,n,\theta)$ in this setting, let us for convenience recall the result from  
the proof of Theorem \ref{abstract}. In particular, we have that 
\begin{equation}\label{norm_sum}
\begin{split}
 \|P_nU^*P_m UP_n - P_nU^*UP_n\| &\leq  \sum_{j=m+1}^\infty \langle  U P_n U^*e_j, e_j\rangle, 
 \qquad m \rightarrow \infty.
 \end{split}
 \end{equation}
 
\begin{theorem}\label{wavelets}
Suppose that $\{\varphi_l\}_{l\in\mathbb{N}}$ is a collection of functions as in (\ref{med}) such that $\mathrm{supp}(\varphi_l) \subset [-T,T]$ for all $l \in \mathbb{N}$ and some $T>0$.
Let $U$ be defined as in Proposition \ref{prop} with $0 < \epsilon \leq \frac{1}{2T}$ and let the bijection
$\rho:\mathbb{N} \rightarrow \epsilon \mathbb{Z}$ defined by 
$\rho(1) = 0, \rho(2) = \epsilon, \rho(3) = -\epsilon, \rho(4) = 2\epsilon, \hdots$. 
For $\theta > 0, n \in \mathbb{N}$ define $\Phi(U,n,\theta)$ 
as in (\ref{thefunction}).  Then, if $\phi, \psi$ satisfy (\ref{wavelet_C}), we have that 
\begin{equation*}
\Phi(U,n,\theta) \leq 
\left(\frac{4\epsilon^{1-2p}\lceil a \rceil C^2}{f(\theta)}\right)^{\frac{1}{2p-1}}
\left(1+\left(\frac{4^p n^{2p}-1}{4^p-1}\right)\right)^{\frac{1}{2p-1}}
= \mathcal{O}\left(n^{\frac{2p}{2p-1}}\right), \qquad n \rightarrow \infty, 
\end{equation*}
where 
$
f(\theta) = (\sqrt{1+4\theta^2} -1)^2/(4\theta^2).$
\end{theorem}
\begin{proof}

To estimate $\Phi(U,n,\theta)$ we will determine bounds on 
$$
\Psi(U,n,\theta) = \min\left\{m \in \mathbb{N}:
\left\|(P_nU^*P_mUP_n\lvert_{P_n\mathcal{H}})^{-1}\right\|\left\|P_nU^*P_mUP_n^{\perp}\right\| \leq \theta\right\}.
$$
Note that if $r < 1$ and 
$
 \|P_nU^*P_m UP_n - P_nU^*UP_n\|  \leq r,
$
then 
$\|(P_nU^*P_mUP_n\lvert_{P_n\mathcal{H}})^{-1}\| \leq \epsilon/(1-\epsilon r)$ (recall that $U^*U = \epsilon^{-1}I$ 
and that $\epsilon \leq 1$). Also, recall (\ref{crep}), so that 
$$
\left\|(P_nU^*P_mUP_n\lvert_{P_n\mathcal{H}})^{-1}\right\|\left\|P_nU^*P_mUP_n^{\perp}\right\| \leq 
\theta
$$ 
when $r$ and $m$ are chosen such that  
$$
\frac{\sqrt{\epsilon r}}{1-\epsilon r} \leq \theta, \qquad \|P_nU^*P_m UP_n - P_nU^*UP_n\|  \leq r,
$$
(recall that $\|U\| = 1/\sqrt{\epsilon}$).
In particular, it follows that 
\begin{equation}\label{jazz}
\Psi(U,n,\theta) \leq \min\{m:    \|P_nU^*P_m UP_n - P_nU^*UP_n\| \leq \epsilon^{-1}(\sqrt{1+4\theta^2} -1)^2/(4\theta^2)   \}.
\end{equation}
To get bounds on $\Psi(U,n,\theta)$ we will proceed as follows. Since $\phi, \psi$ have compact support, it follows that 
$\mathcal{F}\phi, \mathcal{F}\psi$ are bounded. Moreover, by assumption, we have that 
$$
|\mathcal{F}\phi(w)| \leq \frac{C}{|w|^p}, \qquad |\mathcal{F}\psi(w)| \leq \frac{C}{|w|^p}, \qquad \omega \in \mathbb{R}\setminus \{0\}.
$$
And hence, since  
$$
\mathcal{F}\psi_{j,k}(w) = e^{-2\pi \mathrm{i}2^{-j}kw}2^{\frac{-j}{2}}\mathcal{F}\psi(2^{-j}w),
$$
we get that 
\begin{equation}\label{revel}
|\mathcal{F}\psi_{j,k}(w)| \leq 2^{\frac{-j}{2}}\frac{C}{|2^{-j}w|^p}, \qquad \omega \in \mathbb{R}.
\end{equation}
Note that, by the definition of $U$, it follows that 
$$
 \sum_{j=m+1}^\infty \langle  U P_n U^*e_j, e_j\rangle = \sum_{s=m+1}^{\infty}\sum_{t=1}^n 
 |\mathcal{F}\varphi_t(\rho(s))|^2.
$$
And also, by (\ref{revel}) and (\ref{med}) we have,  for $s>0$,
\begin{equation*}
\begin{split}
\sum_{t=1}^n 
 &|\mathcal{F}\varphi_t(\rho(s))|^2 \leq 2\lceil a \rceil|\mathcal{F}\phi(\rho(s)) |^2  + \sum_{j=0}^{\lfloor \log_2(n)\rfloor }
 \sum_{k = -\lceil a \rceil+1}^{2^j\lceil a \rceil-1}|\mathcal{F}\psi_{j,k}(\rho(s))|^2 \\
 &\leq \frac{2\lceil a \rceil C^2}{|\rho(s)|^{2p}} +   \sum_{j=0}^{\lfloor \log_2(n)\rfloor}\sum_{k = -\lceil a \rceil+1}^{2^j\lceil a \rceil-1}2^{-j}\frac{C^2}{|2^{-2j}\rho(s)^2|^p} 
 =  2\lceil a \rceil\left( \frac{C^2}{|\rho(s)|^{2p}} +    \sum_{j=0}^{\lfloor \log_2(n)\rfloor} \frac{C^2}{|2^{-2j}\rho(s)^2|^p}\right) \\
 &\leq \frac{2\lceil a \rceil C^2}{|\rho(s)|^{2p}}\left(1 + \frac{4^pn^{2p} -1}{4^p-1}\right),
\end{split}
\end{equation*}
thus we get that 
\begin{equation}\label{buddy}
\begin{split}
&\sum_{s=m+1}^{\infty}\sum_{t=1}^n |\mathcal{F}\varphi_t(\rho(s))|^2 
\leq 2\lceil a \rceil C^2\left(1 + \frac{4^pn^{2p} -1}{4^p-1}\right)\sum_{s=m+1}^{\infty}\frac{1}{|\rho(s)|^{2p}} \\
 &\leq 2\epsilon^{-2p}2\lceil a \rceil C^2\left(1 + \frac{4^pn^{2p} -1}{4^p-1}\right)\sum_{s=m+1}^{\infty}\frac{1}{s^{2p}} \leq \frac{4\epsilon^{-2p}\lceil a \rceil C^2}{m^{2p-1}}\left(1 + \frac{4^pn^{2p} -1}{4^p-1}\right).
 \end{split}
\end{equation}
Therefore, by using (\ref{norm_sum}) we have just proved that 
$$
\|P_nU^*P_m UP_n - P_nU^*UP_n\| \leq \frac{4\epsilon^{-2p}\lceil a \rceil C^2}{m^{2p-1}}\left(1 + 
\frac{4^pn^{2p} -1}{4^p-1}\right),
$$
and by plugging this bound into (\ref{jazz}) we obtain
$$
\Psi(U,n,\theta) \leq \left(\frac{4\epsilon^{1-2p}\lceil a \rceil C^2}{f(\theta)}\right)^{\frac{1}{2p-1}}
\left(1+\left(\frac{4^p n^{2p}-1}{4^p-1}\right)\right)^{\frac{1}{2p-1}},
$$ which obviously yields the asserted bound on $\Phi(U,n,\theta).$
\end{proof}

The theorem has an obvious corollary for smooth compactly supported wavelets.
\begin{corollary}
Suppose that we have the same setup as in Theorem \ref{wavelets}, and suppose also that 
$\phi, \psi \in C^p(\mathbb{R})$ for some $p \in \mathbb{N}$. Then
$$
\Phi(U,n,\theta) = \mathcal{O}\left(n^{\frac{2p}{2p-1}}\right), \qquad n \rightarrow \infty.
$$
\end{corollary}

\subsection{A Pleasant Surprise}
\begin{figure}
\centering
\includegraphics[height=40mm]{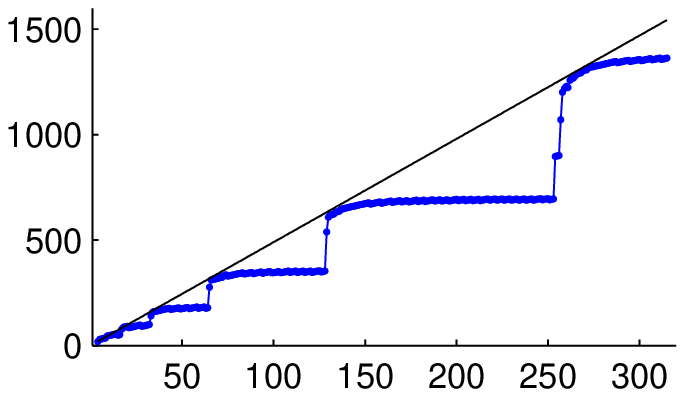}
\includegraphics[height=40mm]{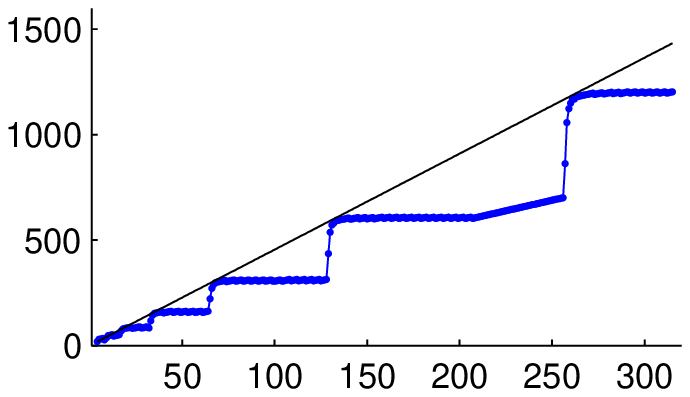}
\caption{The figure shows sections of the graphs of  $\widetilde \Psi (U,\cdot,1)$ (left) and $\widetilde \Psi(U,\cdot,2)$ (right) together 
with the functions (in black) $x \mapsto 4.9x$ (left) and $x \mapsto 4.55x$. In this case $U$ is 
formed by using the Haar wavelets on $[0,1]$.}
\label{fred}
\end{figure}

Note that if $\psi$ is the Haar wavelet and $\phi = \chi_{[0,1]}$ we have that 
$$
|\mathcal{F}\phi(w)| \leq \frac{2}{|w|}, \qquad |\mathcal{F}\psi(w)| \leq \frac{2}{|w|}, \qquad \omega \in \mathbb{R}.
$$
Thus, if we used the Haar wavelets on $[0,1]$ as in Theorem \ref{wavelets} and used the technique in the proof 
of Theorem \ref{wavelets} we would get that 
\begin{equation}\label{quant}
\min\{m:    \|P_nU^*P_m UP_n - P_nU^*UP_n\| =  \epsilon^{-1}(\sqrt{1+4\theta^2} -1)^2/(4\theta^2)   \} = \mathcal{O}\left(n^2\right), \qquad n \rightarrow \infty.
\end{equation}
However, it is tempting to check numerically whether this bound is sharp. Let us denote the quantity 
in (\ref{quant}) by $\widetilde \Psi(U,n,\theta),$ and observe that this can easily be computed numerically. 
Figure \ref{fred} shows $\widetilde \Psi(U,n,\theta)$ for $\theta = 1,2$, where $U$ is defined as in Proposition 
\ref{prop} with $\epsilon = 0.5$. Note that the numerical computation actually shows that 
\begin{equation}\label{on}
\widetilde \Psi(U,n,\theta) = \mathcal{O}\left(n\right),
\end{equation}
which is indeed a very pleasant surprise. In fact, due to the "stair case" growth shown in Figure \ref{fred}, 
the growth is actually better than what 
(\ref{on}) suggests.
The question is whether this is a particular quality of the Haar wavelet, or that one can expect similar behavior of 
other types of wavelets. The answer to this question will be the topic of future papers.

Note that Figure \ref{fred} is interpreted as follows: provided $m \geq 4.9 n$, for example, we can expect this method to reconstruct $g$ to within an error of size $(1+\theta)\| P^{\top}_{n} \beta\|$, where $\theta=1$ in this case.  In other words, the error is only two times greater than the best approximation to $g$ from the finite-dimensional space consisting of the first $n$ Haar wavelets.

\begin{figure}
\centering
\includegraphics[height=40mm]{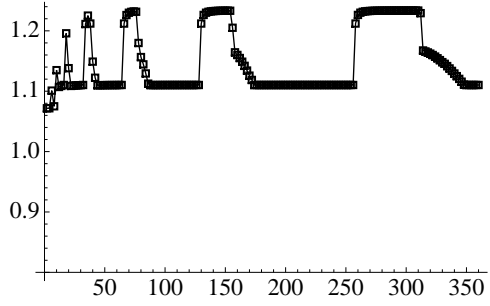}
\caption{The quantity $\| (\epsilon\hat{A})^{-1} \|$ against $n=2,4,\ldots,360$.}
\label{stablefig}
\end{figure}

Having described how to determine conditions which guarantee existence of a reconstruction, in the next section we apply this approach to a number of example problems.  First, however, it is instructive to confirm that these conditions do indeed guarantee stability of the recontruction procedure.  In Figure \ref{stablefig} we plot $\| (\epsilon \hat{A})^{-1} \|$ against $n$ (for $\epsilon = 0.5$), where $\hat{A}$ is formed via (\ref{black}) using Haar wavelets with parameter $m = \lceil 4.9 n \rceil$.  As we observe, the quantity remains bounded, indicating stability. Note the stark contrast to the 
severe instability documented in Figure \ref{freddy33}.

\section{Examples}\label{Numerics}
In this final section, we consider the application of the generalized sampling theorem to several examples.

\subsection{Reconstruction from the Fourier Transform}\label{four}
In this example we consider the following problem.  Let $f \in L^2(\mathbb{R})$ be such that 
$$
f = \mathcal{F}g, \qquad \mathrm{supp}(g) \subset [-T,T].
$$
We assume that we can access point samples of $f$, however, it is not $f$ that is of interest to us, but rather $g$.
This is a common problem in applications, in particular MRI. 
The NS Sampling Theorem assures us that we can recover $g$ from point samples of $f$ as follows:
$$
g = \epsilon \sum_{n = -\infty}^{\infty}f(n\epsilon)\, e^{2\pi \mathrm{i} n\epsilon \cdot }, \qquad \epsilon = \frac{1}{2T},
$$
where the series converges in $L^2$ norm. Note that the speed of convergence depends on how well $g$ can be approximated by the functions $e^{2\pi \mathrm{i} n\epsilon \cdot}$, $n \in \mathbb{Z}.$ 
Suppose now that we consider the function 
$$
g(t) = \cos(2\pi t)\chi_{[0.5,1]}(t).
$$ 
In this case, due to the discontinuity, forming 
\begin{equation}\label{g_N}
g_N = \epsilon \sum_{n = -N}^{N}f(n\epsilon)\, e^{2\pi \mathrm{i} n\epsilon \cdot }, \qquad 
\epsilon = \frac{1}{2}, \quad N \in \mathbb{N},
\end{equation}
may be less than ideal, since the convergence $g_N \rightarrow g$ as $N \rightarrow \infty$ may be slow. 

This is, of course, not an issue if we can access all the samples $\{f(n\epsilon)\}_{n \in \mathbb{Z}},$ however, such 
an assumption may be overly optimistic in application.  Moreover, even if we had access to all samples, we are limited by both processing power and storage to taking only a finite number.

\begin{figure}
\centering
\includegraphics[height=36mm]{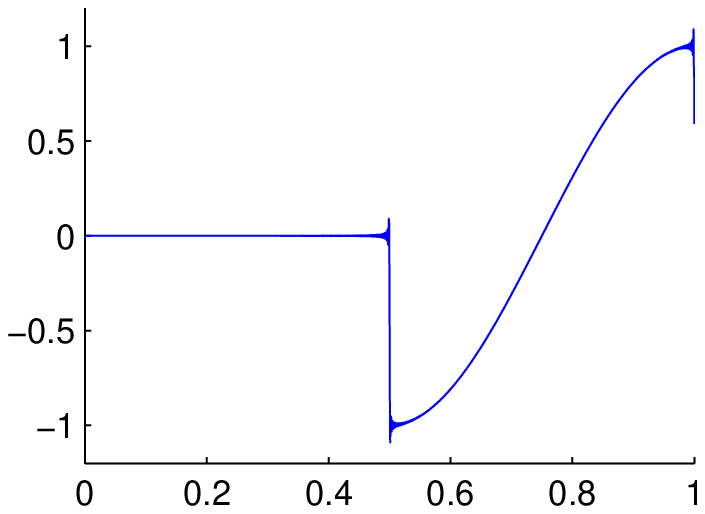}
\includegraphics[height=36mm]{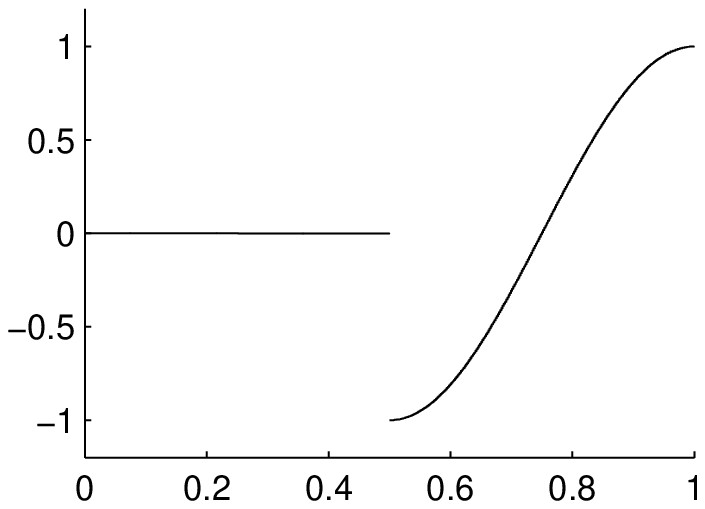}
\includegraphics[height=36mm]{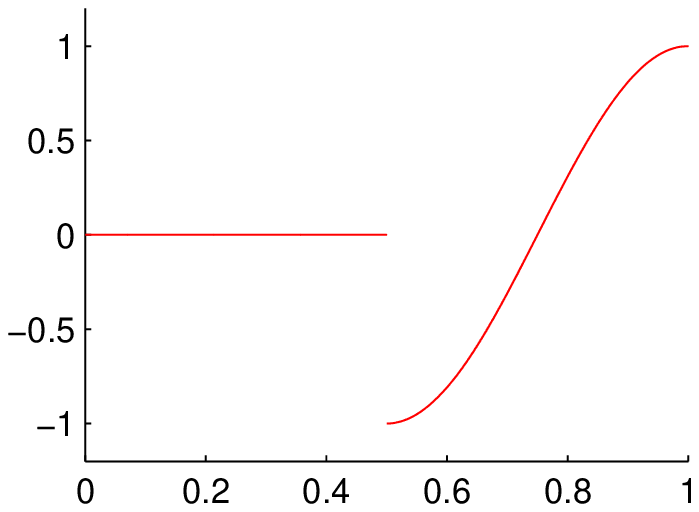}\\
\includegraphics[height=36mm]{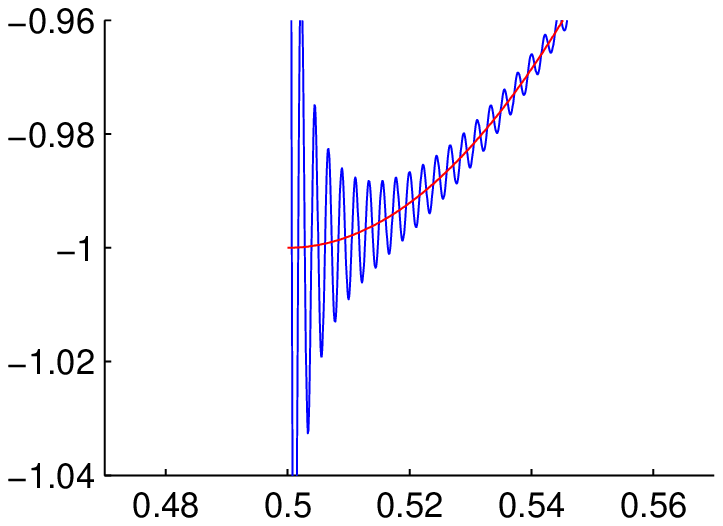}
\includegraphics[height=36mm]{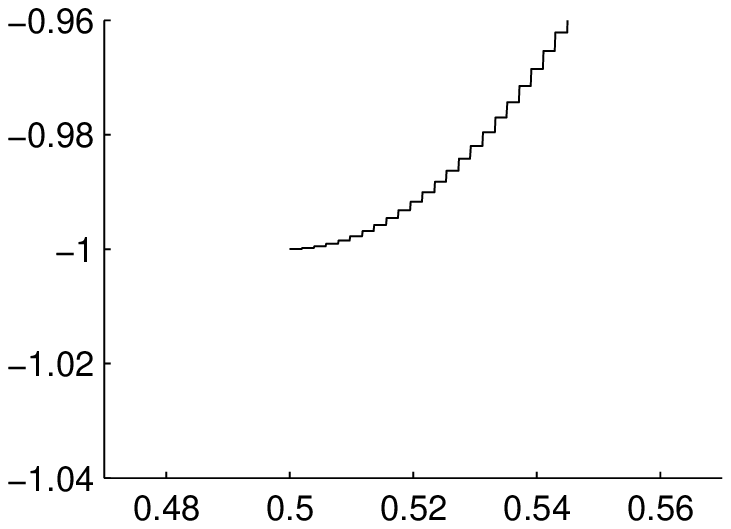}
\includegraphics[height=36mm]{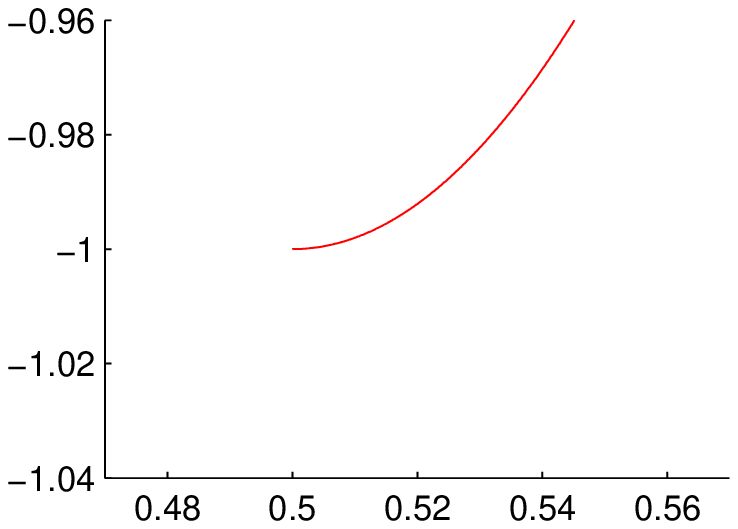}
\caption{The upper figures show $g_N$ (left), $\tilde g_{n,m}$ (middle) and $g$ (right) on the interval $[0,1]$.
The lower figures show $g_N$ (left), $\tilde g_{n,m}$ (middle) and $g$ (right) on the interval $[0.47,0.57]$.}
\label{approx}
\end{figure}

Suppose that we have a more realistic scenario: namely, we are given the finite collection of samples
\begin{equation}\label{eta}
\eta_f = \{f(-N\epsilon),  f((-N+1)\epsilon), \hdots, f((N-1)\epsilon), f(N\epsilon) \},  
\end{equation}
with $N = 900$ and $\epsilon = \frac{1}{2}.$
The task is now as follows: construct the best possible approximation to $g$ based on the vector $\eta_f.$ We can naturally form $g_N$ as in (\ref{g_N}).  This approximation can be visualized in the 
diagrams in Figure \ref{approx}. Note the rather unpleasant Gibbs oscillations that occur, as discussed previously. The problem is simply that 
the set $\{e^{2\pi \mathrm{i} n\epsilon \cdot}\}_{n \in \mathbb{Z}}$ is not a good basis to express $g$ in.
Another basis to use may be the Haar wavelets $\{\psi_j\}$ on $[0,1]$ (we do not claim that this is the optimal basis, 
but at least one that may better capture the discontinuity of $g$).  In particular,  we may express $g$ as
$$
g =  \sum_{j=1}^{\infty} \beta_j \psi_j, \qquad \beta = \{\beta_1, \beta_2, \hdots\} \in l^2(\mathbb{N}).
$$
We will now use the technique suggested in Theorem \ref{Diff_Shan} to construct a better approximation 
to $g$ based on exactly the same input information: namely, $\eta_f$ in (\ref{eta}). Let 
$\widehat U$ be defined as in (\ref{hatU}) with $\epsilon = 1/2$ and let $n = 500$ and $m = 1801.$ In this case 
\begin{equation*}
\begin{split}
\left\|\left(P_n\widehat U^*P_m\widehat UP_n\lvert_{P_n\mathcal{H}}\right)^{-1}\right\| &\leq  0.6169, \\  
\left\|\left(P_n\widehat U^*P_m\widehat UP_n\lvert_{P_n\mathcal{H}}\right)^{-1}P_n \widehat U^*P_m\right\| &\leq 0.7854
\end{split}
\end{equation*}
Define $\tilde \beta = \{\tilde \beta_1, \hdots, \tilde \beta_n\}$ by equation (\ref{black}), and let 
$\tilde g_{n,m} =   \sum_{j=1}^{n}\tilde \beta_j \psi_j.$ The function $\tilde g_{n,m}$ is visualized in Figure \ref{approx}.  Although, the construction of $g_N$ and $\tilde g_{n,m}$ required exactly the same amount of samples of $f$, it is clear from Figure \ref{approx} that $\tilde g_{n,m}$ is favorable.  In particular, approximating $g$ by $\tilde{g}_{n,m}$ gives roughly four digits of accuracy.  Moreover, had both $n$ and $m$ been increased, this value would have decreased.  In contrast, the approximation $g_N$ does not converge uniformly to $g$ on $[0,1]$.

\subsection{Reconstruction from Point Samples}
\begin{figure}
\centering
\includegraphics[height=37mm]{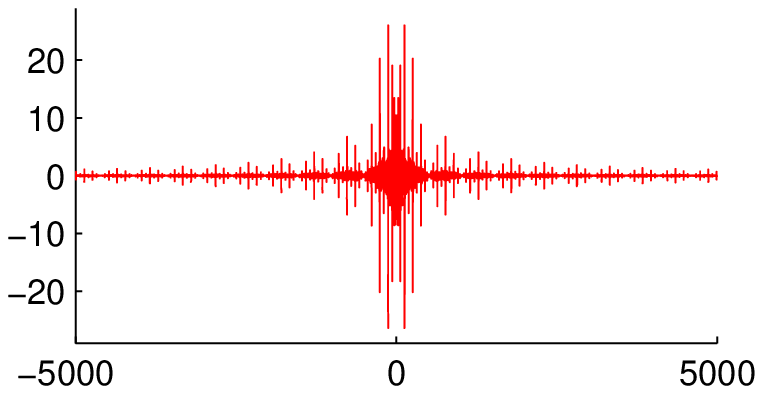}
\includegraphics[height=37mm]{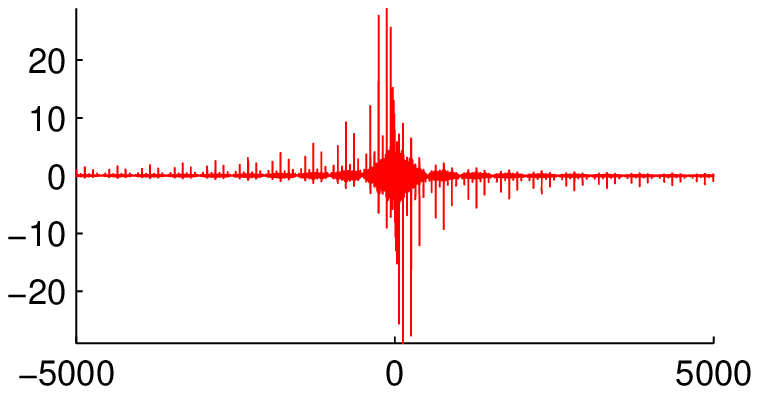}
 \caption{The figure shows $\mathrm{Re}(f)$ (left)  and $\mathrm{Im}(f)$ (right) on the interval $[-5000,5000]$. }
\label{sinc}
\end{figure}

\begin{figure}
\centering
\includegraphics[height=37mm]{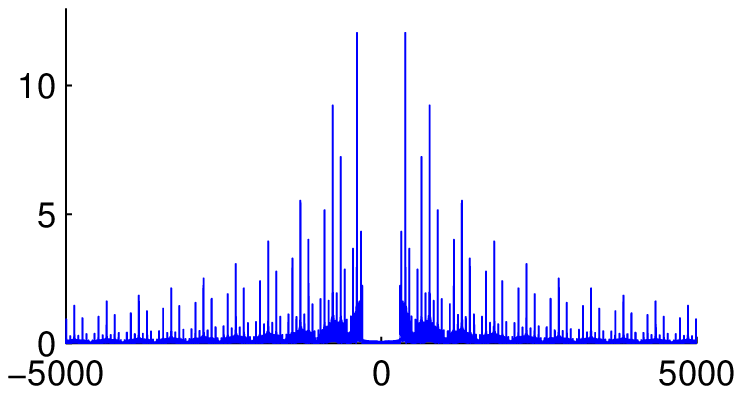}
\includegraphics[height=37mm]{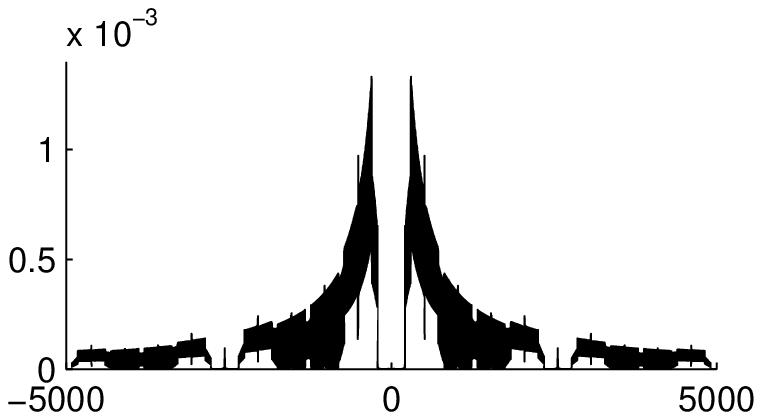}
\caption{The figure shows the error $|f-f_N|$ (left) and $|f - \tilde f|$ (right)
on the interval $[-5000,5000]$.}
\label{sinc2}
\end{figure}

In this example we consider the following problem.  Let 
$f \in L^2(\mathbb{R})$ such that 
$$
f = \mathcal{F}g, \qquad   g(x) 
= \sum_{j=1}^K \alpha_j \psi_j(x) + \sin(2\pi x)\chi_{[0.3,0.6]}(x), 
$$
for $K = 400,$ where $\{\psi_j\}$ are Haar wavelets on $[0,1],$ and 
$\{\alpha_j\}_{j=1}^K$ are some arbitrarily chosen 
real coefficients in $[0,10]$. A section of the graph of $f$ is displayed in 
Figure \ref{sinc}.
The NS Sampling Theorem yields that 
$$
f(t) = \sum_{k = -\infty}^{\infty} f\left(\frac{k}{2}\right) 
\mathrm{sinc}(2t-k),
$$
where the series converges uniformly. Suppose that we can access the 
following pointwise samples of 
$f$:
$$
\eta_f = \{f(-N\epsilon),  f((-N+1)\epsilon), \hdots, 
f((N-1)\epsilon), f(N\epsilon) \}, 
$$
with $\epsilon = \frac{1}{2}$ and $N = 600.$ The task is to reconstruct
 an approximation to $f$ from the samples $\eta_f$ 
in the best possible way. We may of course  form
 $$
 f_N(t) = \sum_{k = -N}^N f\left(\frac{k}{2}\right)
 \mathrm{sinc}(2t-k), 
\qquad N = 600.
 $$
However, as Figure \ref{sinc2} shows, this approximation 
 is clearly less than ideal as $f(t)$ is approximated poorly for large $t$.
It is therefore tempting to try the reconstruction based on Theorem 
\ref{Diff_Shan} and the Haar wavelets on $[0,1]$ (one may of course try a different basis). 
In particular, let 
$$
\tilde f = \sum_{j=1}^n \tilde \beta_j \mathcal{F}\psi_j, \qquad n = 500,
$$ 
where 
$$
\widehat A  \tilde \beta  
=
P_n\widehat U^*P_m
\eta_f, \qquad \widehat A = P_n\widehat U^*P_m\widehat UP_n\lvert_{P_n\mathcal{H}}, 
$$
with $m = 2N+1=1201$
and $\widehat U$ is defined in (\ref{hatU}) with $\epsilon = 1/2$. A section of the errors
$|f - f_N|$ and $|f - \tilde f|$ is shown in Figure \ref{sinc2}. 
In this case we have 
\begin{equation*}
\begin{split}
\left\|\left(P_n\widehat U^*P_m\widehat
    UP_n\lvert_{P_n\mathcal{H}}\right)^{-1}\right\| &\leq  0.9022, \\
\left\|\left(P_n\widehat U^*P_m\widehat
    UP_n\lvert_{P_n\mathcal{H}}\right)^{-1}P_n \widehat U^*P_m\right\|
&\leq 0.9498.
\end{split}
\end{equation*}
In particular, the reconstruction $\tilde f$ is very stable. 
Figure \ref{sinc2} displays how our alternative reconstruction is favorable especially for large $t$.
Note that with the same amount of sampling 
information the improvement is roughly by a factor of ten thousand.

\section{Concluding Remarks}

The framework presented in this paper has been studied via the examples of Haar wavelets and Legendre polynomials.  Whilst the general theory is now well developed, there remain many questions to answer within these examples.  In particular,

\begin{itemize}
\item[(i)] What is the required scaling of $m$ (in comparison to $n$) when the reconstruction basis consists of Legendre polynomials, and how well does the resulting method compare with more well-established approaches for overcoming the Gibbs phenomenon in Fourier series?  Whilst there have been some previous investigations into this particular approach \cite{hrycakIPRM,shizgalGegen2}, we feel that the framework presented in this paper, in particular the estimates proved in Theorem \ref{abstract}, are well suited for understanding this problem.  We are currently investigating this possibility, and will present our results in a future paper. 

\item[(ii)] Whilst Haar wavelets have formed been the principal example in this paper, there is no need to restrict to this case.  Indeed, Theorem \ref{wavelets} provides a first insight into using more sophisticated wavelet bases for reconstruction.  Haar wavelets are extremely simple to work with, however the use of other wavelets presents a number of issues.  In particular, it is first necessary to devise a means to compute the entries of the matrix $U$ in a more general setting.

In addition, within the case of the Haar wavelet, there remains at least one open problem.  The computations in Section \ref{m} suggest that $n \mapsto \Phi(U,n,\theta)$ is bounded by a linear function in this case, meaning that Theorem \ref{wavelets} is overly pessimistic. This must be proven.  Moreover, it remains to be seen whether a similar phenomenon holds for other wavelet bases.

\item[(iii)] The theory in this paper has concentrated on linear reconstruction techniques with full sampling. A natural question is whether one can apply non-linear techniques from compressed sensing to allow for subsampling. Note that, due to the infinite dimensionality of the problems considered here, the standard finite-dimensional techniques are not sufficient. 

\end{itemize}

\section{Acknowledgments}
The authors would like to thank Emmanuel Cand{\`e}s and Hans G. Feichtinger for valuable discussions and input.
\bibliographystyle{abbrv}
\bibliography{bib_file_Shannon}

\end{document}